\documentclass[12pt]{article}
\title{Random Sorting Networks}
\author{Omer Angel \and Alexander E. Holroyd \and Dan Romik \and
B\'alint Vir\'ag}
\date{September 16, 2006}

\usepackage{amsmath,amssymb,amsfonts,amsthm,amstext}
\usepackage{graphicx,color,pst-all,pst-math,sidecap}

\newtheorem{thm}{Theorem}
\newtheorem{conj}{Conjecture}
\newtheorem{prop}[thm]{Proposition}

\newtheorem{eg}[thm]{Example}
\newtheorem{lemma}[thm]{Lemma}

\newtheorem{cor}[thm]{Corollary}

\newcommand{\thmref}[1]{Theorem~\ref{T:#1}}
\newcommand{\lemref}[1]{Lemma~\ref{L:#1}}
\newcommand{\figref}[1]{Figure~\ref{F:#1}}

\newcommand{\ess}{{\cal S}}
\newcommand{\inv}{{\textrm{inv}}}

\newcommand{\sort}{\Omega}

\newcommand{\pusn}{\mathbb{P}_{\textrm{U}}}

\newcommand{\partit}{\textrm{Par}}
\newcommand{\syt}{\textrm{SYT}}
\newcommand{\staircase}{{\triangle}}
\newcommand{\evac}{\Phi}
\newcommand{\EG}{{\textrm{EG}}}

\newcommand{\arch}{\mathfrak{Arch}}
\newcommand{\semi}{\mathfrak{semi}}
\newcommand{\leb}{\mathfrak{Leb}}
\newcommand{\lebb}{\mathcal{L}}

\newcommand{\eps}{\varepsilon}
\newcommand{\cmax}{c_{\text{\rm max}}}
\newcommand{\diff}{\;\Delta\;}
\newcommand{\xlim}{\xrightarrow[n\to\infty]{}}

\renewcommand{\P}{\mathbb{P}}
\renewcommand{\S}{\mathbb{S}}
\newcommand{\E}{\mathbb{E}}
\newcommand{\R}{\mathbb{R}}

\newcommand{\bc}{{\bf c}}
\newcommand{\bu}{{\bf u}}
\newcommand{\bv}{{\bf v}}
\newcommand{\ind}{\mathbf{1}}
\newcommand{\id}{\textrm{id}}
\newcommand{\ntoinf}{\qquad\text{as}\quad n\to\infty}

\newcounter{mycount}
\newenvironment{mylist}{\begin{list}{(\roman{mycount})}
{\usecounter{mycount}\itemsep 0pt}}{\end{list}}

\newcounter{algcount}
\newenvironment{alglist}{\begin{list}{{\bf Step \arabic{algcount}.}\em
}
{\usecounter{algcount}\itemsep 0pt}}{\rm\end{list}}

\begin{document}
\maketitle

\begin{abstract}
A sorting network is a shortest path from $12\cdots n$ to $n\cdots
21$ in the Cayley graph of $S_n$ generated by nearest-neighbour
swaps. We prove that for a uniform random sorting network, as
$n\to\infty$ the space-time process of swaps converges to the
product of semicircle law and Lebesgue measure.  We conjecture
that the trajectories of individual particles converge to random
sine curves, while the permutation matrix at half-time converges
to the projected surface measure of the 2-sphere. We prove that,
in the limit, the trajectories are H\"older-$1/2$ continuous,
while the support of the permutation matrix lies within a certain
octagon. A key tool is a connection with random Young tableaux.
\renewcommand{\thefootnote}{}
\footnote{\hspace{-2em}{\bf Key words:} sorting network, random
sorting, reduced word, maximal chain in the weak Bruhat order,
Young tableau, permutahedron.}
\footnote{\hspace{-2em}{\bf 2000 Mathematics Subject
Classifications:} 60C05; 05E10; 68P10}
\footnote{\hspace{-2em}Funded in part by NSERC discovery grants
(AEH and BV); NSF-FRG grant 0244479 (DR); a Connaught grant, a
Canada Research Chair and a Sloan Fellowship (BV); a PIMS
postdoctoral fellowship (OA); and by MSRI, PIMS and BIRS (all
authors).}
\end{abstract}

\section{Introduction}

\begin{figure}
\begin{center}
  \includegraphics[width=\textwidth]{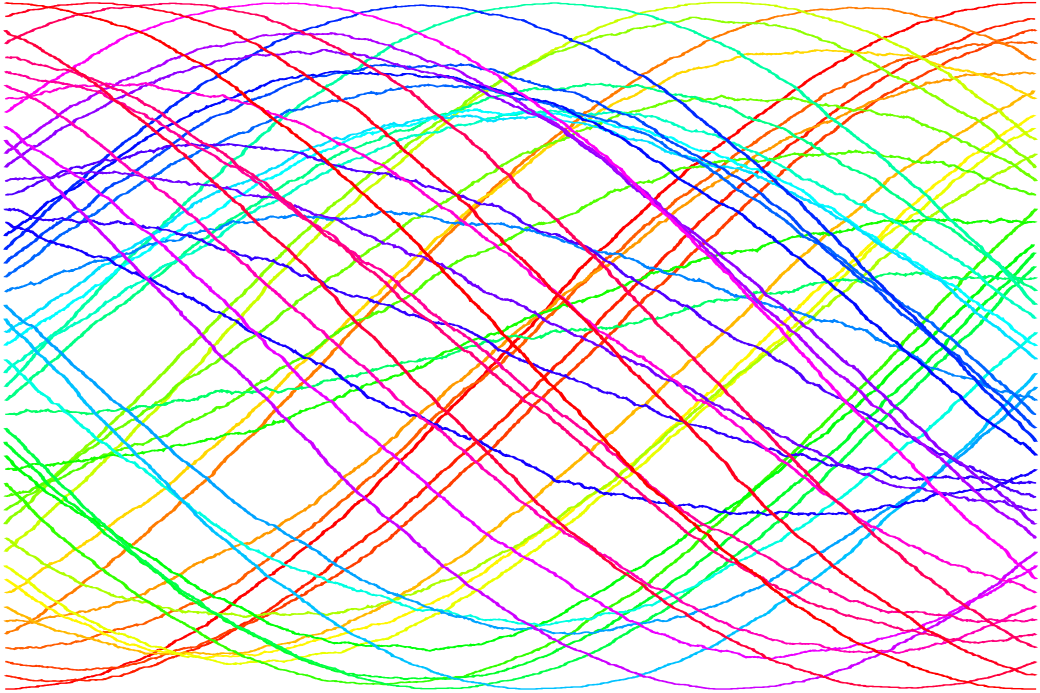}
  \caption{Selected particle trajectories for a uniformly chosen
    2000-element sorting network.}
  \label{F:traj}
\end{center}
\end{figure}
\begin{SCfigure}
  \includegraphics[height=7cm]{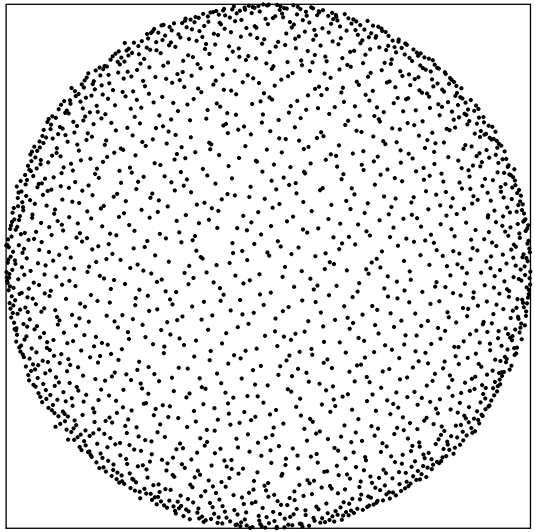}
  \caption{The permutation
  matrix of the half-time configuration $\sigma_{N/2}$ for a uniformly
  chosen 2000-element sorting network.}
  \label{F:half}
\end{SCfigure}

\begin{figure}[t]
\begin{center}
  \includegraphics[height=3.5cm]{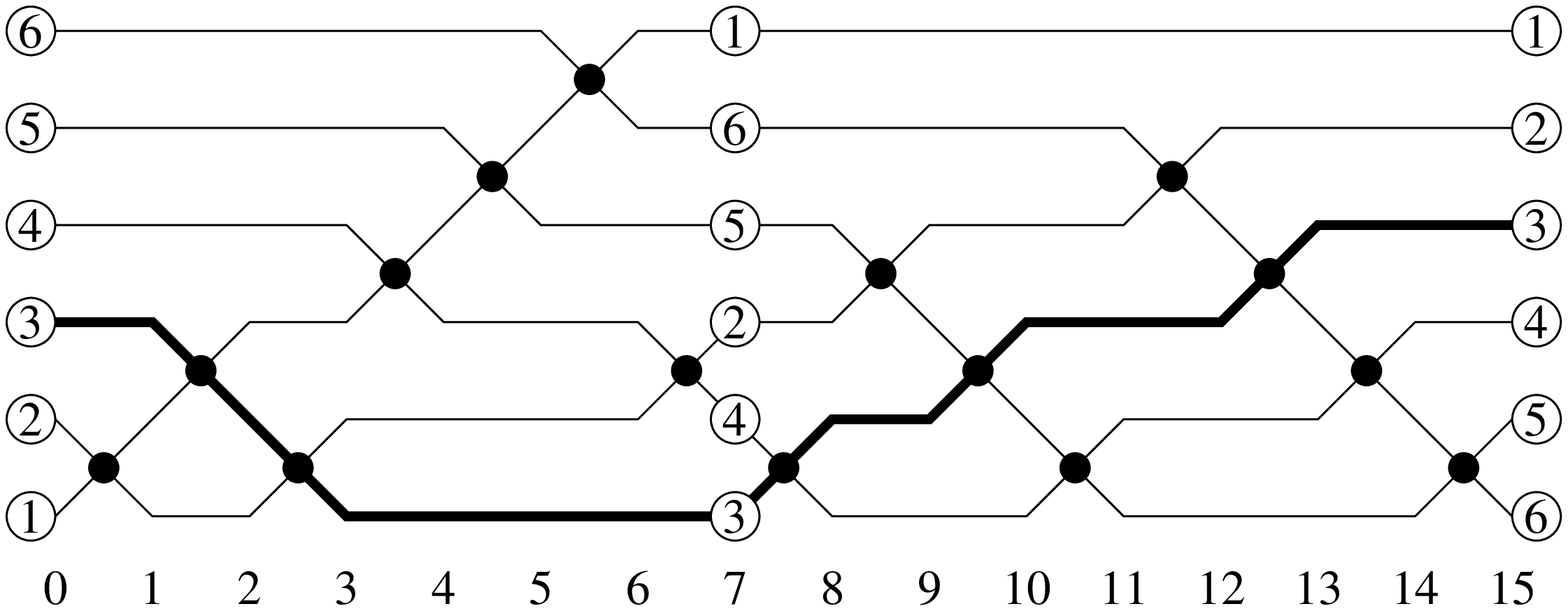}
  \hfill
  \includegraphics[height=3.5cm]{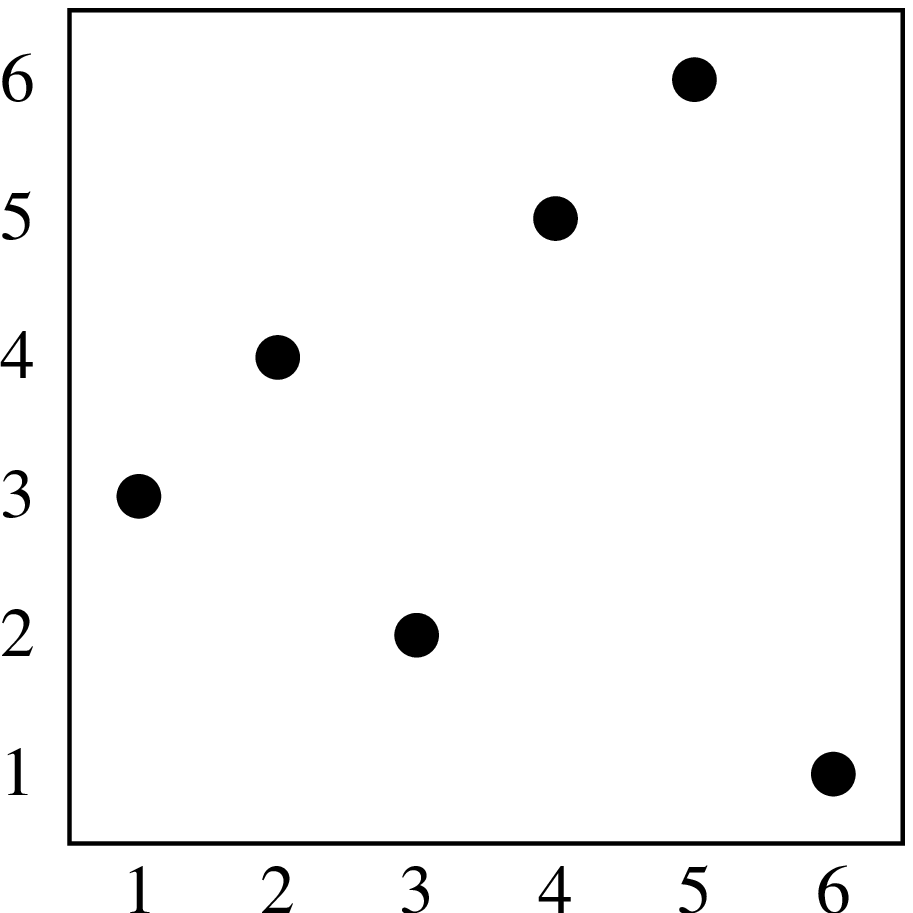}
  \caption{{\em Left:} the ``wiring diagram'' of the 6-element sorting
    network $\omega=(1,2,1,3,4,5,2,1,3,2,1,4,3,2,1)$. The swap process is
    shown by the black discs. The trajectory of particle 3 is highlighted.
    {\em Right:} the graph (or permutation matrix) of the configuration
    $\sigma_7=(3,4,2,5,6,1)$ of $\omega$ at time 7.}
  \label{F:wire}
\end{center}
\end{figure}

\enlargethispage*{1cm} Let $\ess_n$ be
the symmetric group of all permutations
$\sigma=(\sigma(1),\ldots,\sigma(n))$ on $\{1,\ldots,n\}$, with
composition given by $(\sigma\tau)(i):=\sigma(\tau(i))$. For
$1\leq s \leq n-1$ denote the adjacent transposition or {\bf swap}
at location $s$ by $\tau_s := (\,s\ \
s+1\,)=(1,2,\ldots,s+1,s,\ldots,n)\in\ess_n$. Denote the {\bf
identity} $\id:=(1,2,\ldots,n)$ and the {\bf reverse permutation}
$\rho:=(n,\ldots,2,1)$. An $n$-element {\bf sorting network} is a
sequence $\omega=(s_1,\ldots,s_N)$ such that
\[
\tau_{s_1} \tau_{s_2} \cdots \tau_{s_N}=\rho
\]
where
\[
N:=\binom{n}{2}.
\]
(It is easily verified that $N$ is the minimum possible length of
a sequence of swaps whose composition is $\rho$, while $\rho$ is
the unique permutation for which this minimum length is
maximized.) For $1\leq k\leq N$ we refer to $s_k=s_k(\omega)$ as
the $k$th {\bf swap location}, and we call the permutation
$\sigma_k=\sigma_k(\omega):=\tau_{s_1} \cdots \tau_{s_k}$ the {\bf
configuration at time $k$}. We call $\sigma_k^{-1}(i)$ the {\bf
location of particle $i$} at time $k$, and we call the function
$k\mapsto \sigma_k^{-1}(i)$ the {\bf trajectory} of particle $i$.
See Figures \ref{F:traj}, \ref{F:half} and \ref{F:wire} for some
illustrations.

Let $\sort_n$ be the set of all $n$-element sorting networks, and
let $\pusn=\pusn^n$ denote the uniform probability measure on
$\sort_n$ (assigning probability $1/\#\sort_n$ to each
$\omega\in\sort_n$).  We refer to a random sorting network chosen
according to $\pusn$ as a {\bf uniform sorting network (USN)}.

Our first results concern the swap locations.
\enlargethispage*{1cm}
\begin{thm}[Stationarity and semicircle law]
 \label{T:stat-semi}
 Let $\omega_n$ be a uniform $n$-element sorting network.
\begin{mylist}

\item The random sequence $(s_1,\ldots,s_N)$ of swap locations is
stationary; that is $(s_1,\ldots,s_{N-1})$ and $(s_2,\ldots ,s_N)$
are equal in law under $\pusn$.

\item The first swap location $s_1$ satisfies the
convergence in distribution
\[
2 s_1(\omega_n)/n - 1 \Longrightarrow Z \ntoinf
\]
 where $Z$ is a random variable with semicircle law;
that is with probability density function $\frac{2}{\pi}
\sqrt{1-y^2}$ for $y\in(-1,1)$.
\end{mylist}
\end{thm}

\begin{figure}
\begin{center}
  \includegraphics[height=7.5cm]{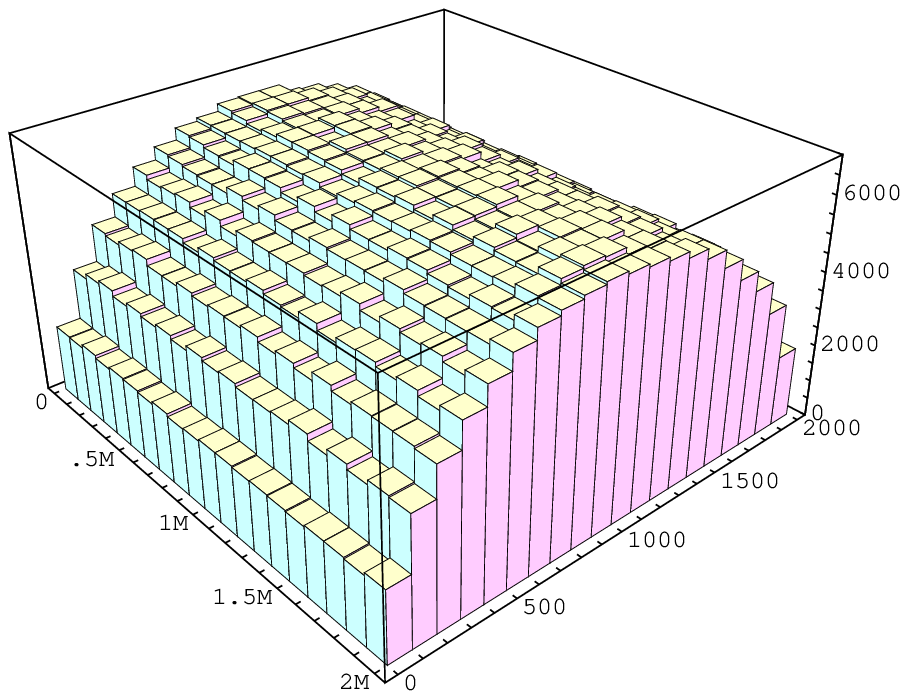}
  \caption{A histogram of the swap process for a uniformly chosen
    2000-element sorting network. (The height of each column represents the
    number of swaps in the corresponding space-time
    window.)}
  \label{F:eta}
\end{center}
\end{figure}

In fact we can compute the exact distribution of $s_1$ for each
$n$; see Proposition \ref{semi}.  In addition we establish the
following ``law of large numbers'' for the swap locations. For an
$n$-element sorting network $\omega$, define the {\bf scaled swap
process} $\eta=\eta(\omega)$ to be the measure
\[
\eta := \frac 1N \sum_{k=1}^N \delta
                              \Big(\frac kN,\; \frac {2 s_k}{n}-1\Big),
\]
where $\delta(x,y)$ is the point measure at $(x,y)$ on ${\mathbb
R}^2$. \figref{eta} is a histogram of $\eta$ for a uniform
2000-element sorting network.  Denote the semicircle measure by
$\semi(dy):=\frac{2}{\pi}\sqrt{1-y^2} \ind_{y\in(-1,1)}\, dy$, and
Lebesgue measure on $[0,1]$ by $\leb(dx):=\ind_{x\in[0,1]}\,dx$.

\enlargethispage*{1cm}
\begin{thm}[Law of large numbers] \label{T:LLN}
Let $\omega_n$ be a uniform $n$-element sorting network.  The scaled
swap process $\eta$ satisfies
\[
\eta(\omega_n) \Longrightarrow \leb \times \semi \ntoinf.
\]
 Here $\Longrightarrow$ denotes convergence in distribution of random
measures in the vague topology on Borel measures on $\R^2$, and the right
side denotes the deterministic product measure.
\end{thm}

For a sorting network $\omega$, define the {\bf scaled trajectory}
$T_i(t)=T_i(t,\omega)$ of particle $i$ by
\[
T_i(t):=2\sigma_{ tN}^{-1}(i)/n-1
\]
when $tN$ is an integer, and by linear interpolation for other
$t\in[0,1]$.

\begin{thm}[H\"older trajectories] \label{T:Holder}
  Let $\omega_n$ be a uniform $n$-element sorting network.
  \begin{mylist}
  \item For any $\varepsilon>0$, the scaled trajectories satisfy
    \begin{equation*}
      \pusn^n
      \Big( \forall i,s,t:\;
      \big| T_i(t)-T_i(s) \big|  \le  \sqrt{8} \big| t-s \big|^{1/2} + \eps
      \Big)
      \to 1 \ntoinf.
    \end{equation*}
  \item
    Let $T_{i(n)}$ be the scaled trajectory of an arbitrarily chosen
    particle $i(n)=i(n,\omega_n)$ in $\omega_n$. Then the random sequence
    $\{T_{i(n)}\}_{n=1}^\infty$ has subsequential limits in distribution
    with respect to uniform convergence of functions, and any subsequential
    limit is supported on $\mbox{\rm H\"older}(\sqrt{8},\tfrac12)$ continuous paths.
  \end{mylist}
\end{thm}

\begin{figure}
\begin{center}
  \includegraphics[width=\textwidth]{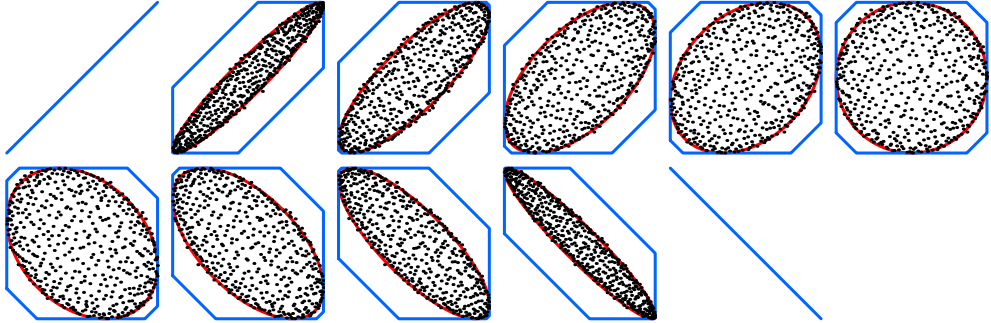}
  \caption{Graphs of the configurations at times $0,\tfrac N{10},
    \tfrac{2N}{10} ,\ldots, N$ for a uniformly chosen 500-element sorting
    network. Also shown are the asymptotic ``octagon bounds'' of
    \thmref{octagon}, and the conjectural asymptotic ``ellipse bounds''
    implied by Conjecture \ref{arch-conj}.}
  \label{fig-oct}
\end{center}
\end{figure}

For a uniform sorting network, the particle configuration
$\sigma_k$ at a given time is a random permutation.  We prove the
following bounds on its distribution.
\begin{thm}[Octagon bounds]\label{T:octagon}
Let $\omega_n$ be a uniform $n$-element sorting network. For any
$\eps>0$ we have
\[
\pusn^n\bigg( \forall k,i:\;
\begin{aligned}
    \big|\sigma_k(i)-i\big|     &< d_k+\eps n, \quad\text{\rm and}\;   \\
    \big|\sigma_k(i)-(n-i)\big| &< d_{N-k}+\eps n
\end{aligned}
\bigg)\to 1\ntoinf,
\]
\\[-2mm]
where $d_k := n \sqrt{\tfrac{k}{N}\left(2-\tfrac{k}{N}\right)}$.
\end{thm}
\thmref{octagon} states that for each $t$, all the 1's in the
permutation matrix of the configuration $\sigma_{\lfloor
tN\rfloor}$ lie within a certain octagon asymptotically almost
surely; see Figure \ref{fig-oct}.

Results of \cite{edelmangreene} and \cite{greeneetal} give rise to
an efficient algorithm for exactly sampling a uniform sorting
network (specifically, see Theorems \ref{T:EG_bij} and
\ref{hookwalk} in this article). The resulting simulations,
together with heuristic arguments, have led us to striking
conjectures about the asymptotic behaviour of the uniform sorting
network.

\pagebreak
\figref{traj} illustrates some trajectories for a
uniform 2000-element sorting network.  We conjecture that as
$n\to\infty$, all particle trajectories converge to sine curves of
random amplitudes and phases.

\begin{conj}[Sine trajectories]
\label{sine-conj}
 Let $\omega_n$ be an $n$-element uniform sorting
network and let $T_i$ be the scaled trajectory of particle $i$.
For each $n$ there exist random variables
$(A^n_i)_{i=1}^n,(\Theta^n_i)_{i=1}^n$ such that for all
$\eps>0$,
\[
  \pusn^n \Big( \max_{i\in[1,n]} \max_{t\in[0,1]}
    \big| T_i(t,\omega_n)- A^n_i \sin(\pi t+\Theta^n_i) \big| > \eps
  \Big) \to 0 \ntoinf.
\]
\end{conj}

Figures \ref{F:half} and \ref{fig-oct} illustrate the graphs
$\{(i,\sigma_k(i) ):i\in[1,n]\}$ (i.e. the locations of 1's in the
permutation matrix) of some configurations from uniform sorting
networks. We conjecture that as $n\to\infty$ the graphs
asymptotically concentrate in a family of ellipses, with a certain
particle density in the interior of the ellipse.
 Define the {\bf scaled configuration} $\mu_t=\mu_t(\omega)$
at time $t$ by
\begin{equation} \label{eq:mutdef}
  \mu_t:=\frac 1n \sum_{i=1}^n \delta\Big(\frac {2i}{n}-1\;,\;\frac {
    2\sigma_{\lfloor tN\rfloor}(i)}{n}-1\Big).
\end{equation}
We define the {\bf Archimedes measure} with parameter $t\in(0,1)$
by
$$\arch_t(dx\times dy):=\frac{1}{2\pi} \;
\sqrt{ \Big[{\sin^2(\pi t) +2xy \cos(\pi t) -x^2-y^2}\Big]^{-1}
\vee 0}\;\;dx\;dy.$$

\pagebreak
\begin{conj}[Archimedes configurations]
\label{arch-conj}
Let $\omega_n$ be an $n$-element uniform sorting network.
For all $t\in(0,1)$, the scaled configuration at time $t$ satisfies
$$\mu_t(\omega_n) \Longrightarrow \arch_t
\ntoinf.$$ Here
$\Longrightarrow$ denotes convergence in distribution in the vague
topology for random Borel measures on $\R^2$.
\end{conj}

In the case $t=1/2$, the measure $\arch_{1/2}$ has density
$1/\big(2\pi\sqrt{1-x^2-y^2}\big)$ on the circular disc
$x^2+y^2<1$. This is the unique circularly symmetric measure whose
linear projections are uniform. It may be obtained by projecting
surface area measure on the 2-sphere in $\R^3$ onto $\R^2$---that
this gives a measure with the aforementioned property follows from
the observation of Archimedes that the surface area of a sphere
between two horizontal planes equals the corresponding area of a
circumscribed vertical cylinder.  (The claimed uniqueness follows
from uniqueness of the characteristic function, \cite[Theorem
5.3]{kallenberg}). For general $t$ the measure $\arch_t$ is
obtained from $\arch_{1/2}$ by the linear transformation $(x,y)
\mapsto (x,x \cos(\pi t)+y \sin(\pi t))$, and is supported on the
interior of an ellipse---see Figure~\ref{fig-oct}.

\begin{figure}
\begin{center}
  \includegraphics[height=6cm]{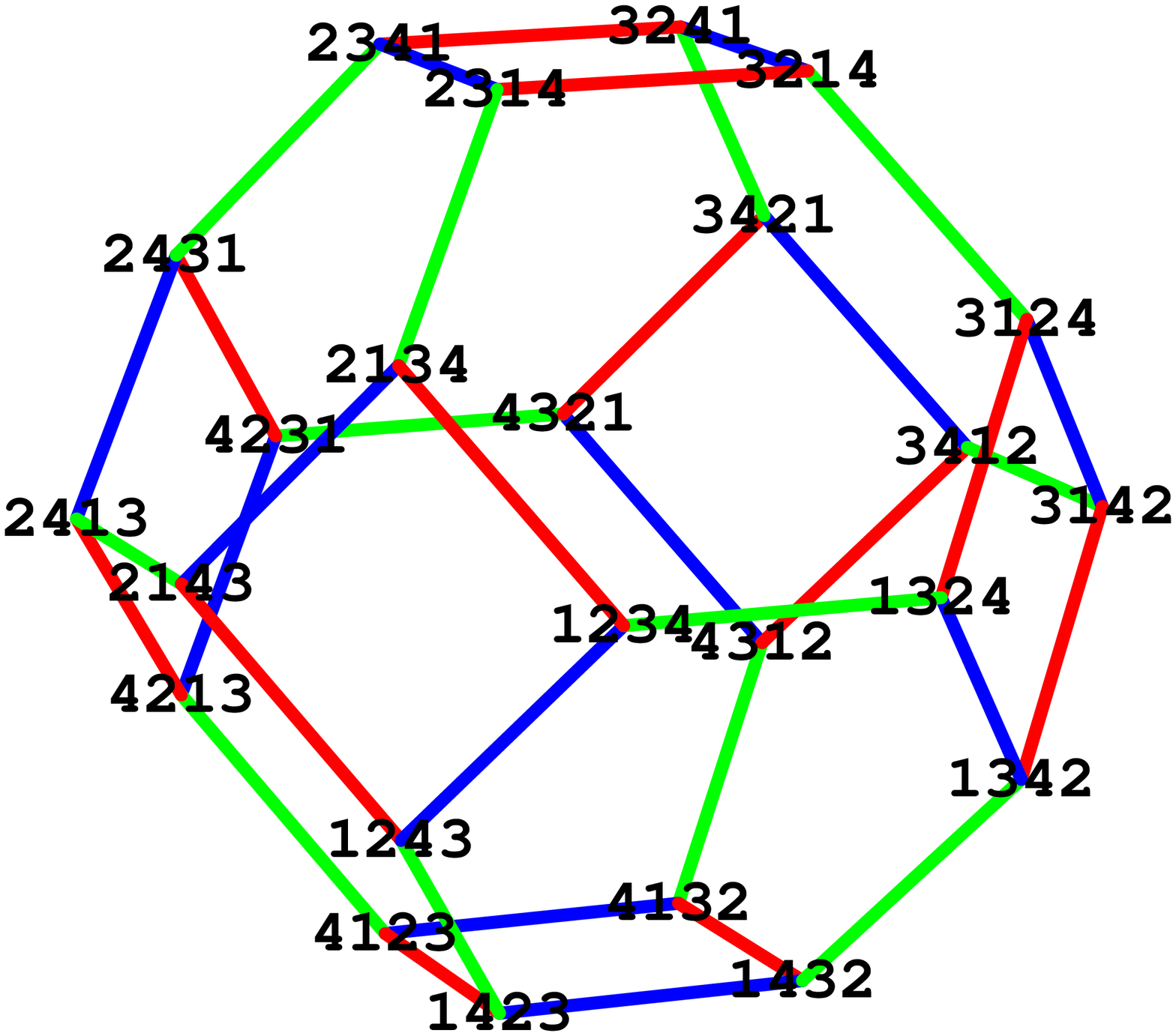}
  \caption{The permutahedron for $n=4$.}
  \label{fig-perm}
\end{center}
\end{figure}

Conjectures \ref{sine-conj} and \ref{arch-conj} (and more) are
implied by a very natural conjecture about the geometry of uniform
sorting networks.  The {\bf permutahedron} is the natural
embedding of the Cayley graph $({\cal S}_n, (\tau_i)_{i=1}^{n-1})$
in Euclidean space in which we assign the permutation
$\sigma\in{\cal S}_n$ to the point
$$\sigma^{-1}=(\sigma^{-1}(1),\ldots,\sigma^{-1}(n))\in\R^n.$$
For all $\sigma\in{\cal S}_n$, clearly $\sigma^{-1}$ lies on the
$(n-2)$-sphere
$${\mathbb S}_n:=\Big\{
{\textstyle z\in \R^n: \sum_{i=1}^n z_i=\tfrac{n(n+1)}{2} }\Big\}
\cap \Big\{{\textstyle z\in \R^n :\sum_{i=1}^n
z_i^2=\tfrac{n(n+1)(2n+1)}{6}}\Big\},$$ while $\id^{-1}$ and
$\rho^{-1}$ are antipodal points on ${\mathbb S}_n$. Furthermore
each edge of the Cayley graph has Euclidean length
$\|\sigma^{-1}-(\sigma\tau_i)^{-1}\|_2=\surd 2$.  See Figure
\ref{fig-perm} for an illustration of the case $n=4$ (where
${\mathbb S}_4$ is a 2-sphere).  A sorting network corresponds to
a shortest path from $\id^{-1}$ to $\rho^{-1}$ in the Cayley
graph.  It is natural to guess that such a path might typically be
close to a {\bf great circle} of ${\mathbb S}_n$; that is, a
Euclidean circle in $\R^n$ having the same centre and radius as
${\mathbb S}_n$.  We show that, if a sorting network lies close to
some great circle, then its trajectories are approximately sine
curves, its particle configurations approximate Archimedes
measure, and its swap locations are approximately governed by the
semicircle law.

\begin{thm}[Great circles] \label{T:great-circle}
Suppose for each $n$ that $\omega_n$ is a (non-random)
$n$-element sorting network, and suppose that there is a sequence
of great circles $c_n\subset {\mathbb S}_n$ such that
$$d_\infty(\omega_n,c_n)=o(n) \quad\text{as}\quad n\to\infty,$$
(with distance defined as $d_\infty(\omega,c):=
\max_{i\in[1,N]}\inf_{z\in c} \|\sigma^{-1}_i-z\|_\infty$). Then:
\begin{mylist}
\item there exist $a_i^n,\theta_i^n$ such that the scaled trajectories satisfy
$$\max_{i\in[1,N]}\max_{t\in[0,1]}\Big| T_i(t,\omega_n)- a^n_i \sin(\pi
t+\theta^n_i)\Big| \to 0 \ntoinf;$$
\item for all $t\in(0,1)$, the scaled configuration satisfies the
vague convergence
$$\mu_t(\omega_n) \Longrightarrow \arch_t
\ntoinf;$$
\item the scaled swap process satisfies the vague convergence
$$\eta(\omega_n) \Longrightarrow \leb \times \semi
\ntoinf.$$
\end{mylist}
\end{thm}

We conjecture that, asymptotically almost surely as $n\to\infty$,
the uniform sorting network does indeed lie close to a great circle on
the permutahedron.
\begin{conj}[Great circles]\label{conj-great-circle}
Let $\omega_n$ be an $n$-element uniform sorting network. For each
$n$ there exists a random great circle $C_n\subset {\mathbb S}_n$
such that
$$d_\infty(\omega_n,C_n)=o(n) \quad\text{in probability as }n\to\infty.$$
\end{conj}

Simulations provide overwhelming numerical evidence in support of
Conjecture \ref{conj-great-circle}. Indeed, the evidence suggests
that for the optimum great circle, typically
$d_\infty(\omega_n,C_n)\approx \mbox{const} \times n^\alpha$,
where $\alpha\approx 1/2$. For example, an exact simulation of a
10000-element uniform sorting network gave
$d_\infty(\omega_{10000},c)\leq 159$ for a certain great circle
$c$. If Conjecture \ref{conj-great-circle} holds then, by
\thmref{great-circle}, Conjectures \ref{sine-conj} and
\ref{arch-conj} follow, as well as the result in \thmref{LLN}. The
fact that \thmref{LLN} does indeed hold thus provides some further
circumstantial evidence for Conjecture \ref{conj-great-circle}. It
is interesting that the proofs of \thmref{LLN} and
\thmref{great-circle}(iii) use entirely different methods.

In addition to \thmref{octagon}, we note that certain other special
permutations may be shown to have asymptotically much lower probability
than others. Since the number of permutations is $n!$, at any given time
step $k\in[0,N]$ there must exist some permutation which is visited with
probability at least $1/n! \ge \exp[-n\log n]$. However, some permutations
are much less likely, as illustrated by the following.
\begin{eg}
\label{double-flip} For $n$ even, let  $h=N/2-n/4$, and consider
the permutation $\psi:=\Big(\tfrac n2,\tfrac n2
-1,\ldots,1,\;\;n,n-1,\ldots,\tfrac n2 +1\Big).$
The probability that the uniform sorting network passes through
particle configuration $\psi$ equals
\[
\pusn^n(\sigma_h =\psi)= \exp\big[-\tfrac{\log 2}{4}n^2 + O(n)\big]
\ntoinf.
\]
(This will be verified in Section \ref{Young tableaux}.)
\end{eg}

\subsection*{Remarks}

\sloppy
\paragraph{History and connections.} \ Sorting networks were first
considered by Stanley \cite{stanleypaper}, who proved the
remarkable formula \fussy
\begin{equation} \label{eq:numsortnet}
\#\sort_n = \frac{\binom{n}{2}!}{1^n 3^{n-1} 5^{n-2} \cdots
(2n-3)^1}.
\end{equation}
Another breakthrough was achieved by Edelman and Greene
\cite{edelmangreene}, who obtained a {\em bijective} proof of
(\ref{eq:numsortnet}). (A related approach to the enumeration of sorting
networks was independently developed by Lascoux and Sch\"utzenberger; see
\cite{lascouxschutz}, \cite[p.\ 94--95]{garsia}.) The Edelman-Greene
bijection is between the set $\sort_n$ of sorting networks and the set of
all staircase-shape standard Young tableaux of size $n$. This bijection
will be an important ingredient for our results; we describe it in Section
\ref{sec:eg}. See \cite{felsner, garsia,goodman, little, reiner} for
further background.

Sorting networks are of interest in computer science, since they
can be interpreted as networks of comparators capable of sorting
any sequence into descending order; see \cite[Exer.\
5.3.4.36--38]{knuth}. There is also a connection with
change-ringing (English-style church bell ringing); for background
see \cite{white} and the references therein.

\paragraph{About the proofs.} The proof of \thmref{stat-semi}(i)
is very simple, and the proof of (ii) is straightforward given the
results of \cite{edelmangreene}. Similar computations appear in
\cite{reiner}. Our proofs of Theorems~\ref{T:LLN}, \ref{T:Holder}
and \ref{T:octagon} are more involved, and depend on results from
\cite{pittelromik} on limiting profiles for random Young tableaux.
A key tool is an extension of the result in \cite{pittelromik}
from square tableaux to staircase tableaux; see Section
\ref{sec:limit}.  \thmref{Holder} is a straightforward consequence
of \thmref{octagon} together with \thmref{stat-semi}(i).  The
proof of \thmref{great-circle} employs geometric arguments, and
relies on the characterization of $\arch_{1/2}$ as the unique
measure all of whose linear projections are uniform on $[-1,1]$.

\paragraph{Simulations.} As remarked above, simulation evidence
strongly supports Conjecture \ref{conj-great-circle}. The measurement
$d_\infty(\omega_{10000},c)\leq 159$ was obtained by using an exact
simulation of a 10000-element USN $\omega_{10000}$, and calculating the
maximum $L^\infty$ distance from the configuration $\sigma_k^{-1}$ at time
$k$ to a point moving at constant angular speed around the great circle $c$
that passes through $\sigma_0^{-1}$ and $\sigma_{N/2}^{-1}$. In contrast,
applying the same procedure to the ``bubble sort'' network
$\omega=(1,2,\ldots,n,\; 1,2,\ldots,n-1,\; \ldots,\; 1,2,\; 1)$ gives for
$n=10000$ a distance of approximately $9997$. It is also easy to see that
the condition $d(\omega_n,c_n)=o(n)$ does not hold for every sequence of
sorting networks $\omega_n$. For example, it does not hold for any sequence
of sorting networks which pass through the permutation $\psi$ in Example
\ref{double-flip}, since the configuration at time $\lfloor N/2\rfloor
\approx h$ cannot satisfy the condition in \thmref{great-circle}(ii).

\begin{SCfigure}
  \includegraphics[height=7cm]{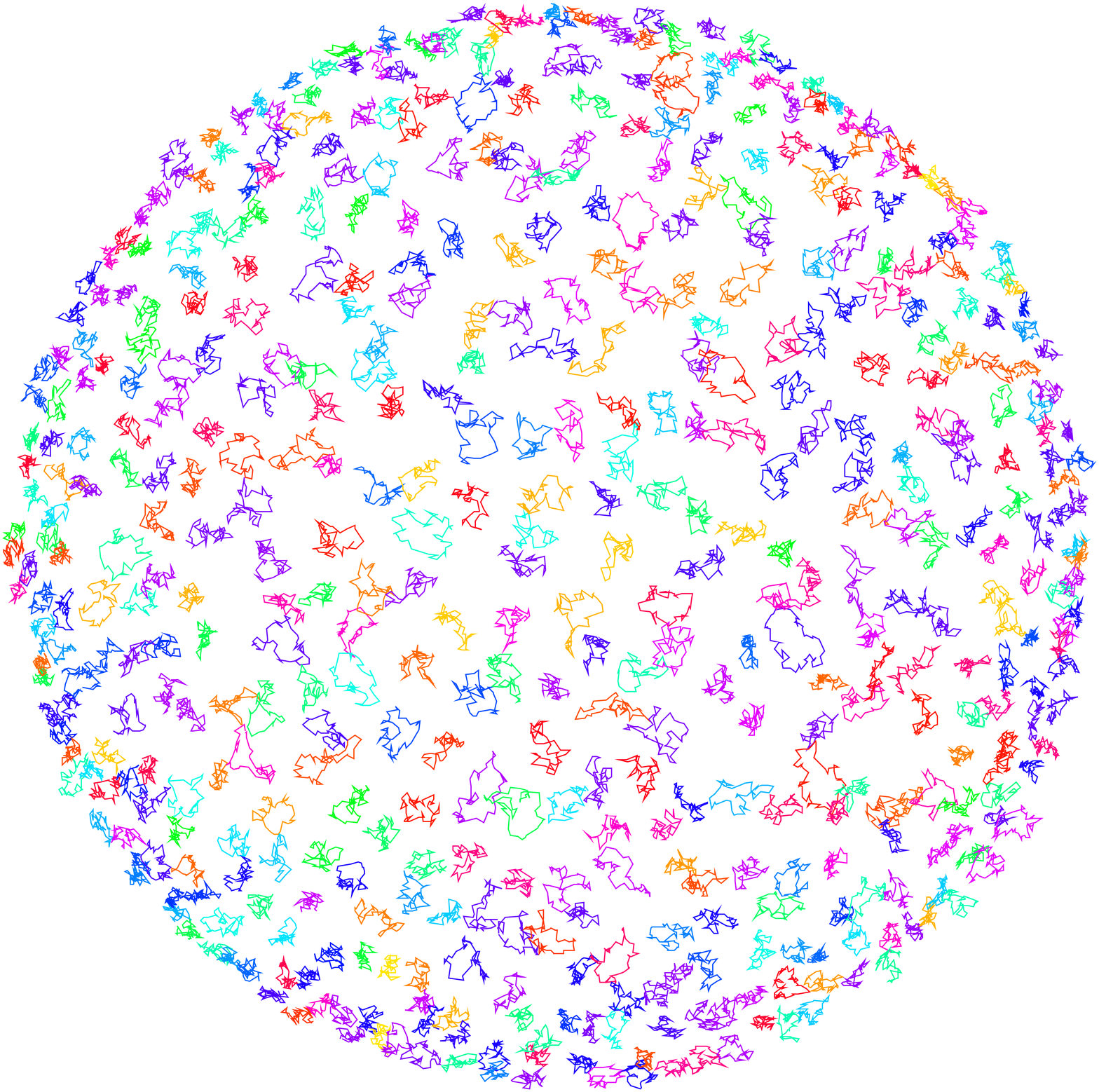}
  \caption{The evolution of the permutation graph of a sliding window, modulo
  uniform rotation, for a
  uniformly chosen
    500-element sorting network.}
  \label{sliding}
\end{SCfigure}

A particularly striking illustration of Conjectures
\ref{sine-conj}--\ref{conj-great-circle} results from plotting the graph of
the permutation $\sigma_{k}^{-1} \sigma_{k+N/2}^{\ }$, and then viewing the
animation as $k$ varies. Stationarity (\thmref{stat-semi}(i)) implies that
at any given time the picture will resemble Figure \ref{F:half}, while at
time $k=N/2$ the initial picture will have been exactly rotated by $\pi/2$.
In fact (for large $n$) the points appear to rotate all at the same
constant angular speed. To further illustrate this we may simultaneously
rotate the entire picture by the (uniformly changing) angle $-\pi k/N$, and
plot the resulting paths of the moving points as $k$ increases from $0$ to
$N/2$. This is shown in Figure~\ref{sliding}. The observation that each
path is localized is a manifestation of Conjectures \ref{sine-conj} and
\ref{conj-great-circle}.

\paragraph{Further works.} In forthcoming articles
  \cite{sort-pen,sort-local,sort-weak} we study several closely related
issues. In \cite{sort-weak} we prove further bounds on the configurations
$\sigma_k$ in the USN. In \cite{sort-local} we study the {\em local}
structure of the swap process. In \cite{sort-pen} we study another natural
probability measure on sorting networks, in which at every step, a swap
location is chosen uniformly from among those locations where the two
particles are in increasing order. It turns out that this model can be
analyzed in detail via the theory of exclusion processes. Its behaviour is
very different from that of the USN, but it has the property, apparently
shared by the USN (see Conjecture \ref{sine-conj}), that asymptotically
each particle initially moves at a well-defined randomly chosen speed, and
continues on a trajectory which is deterministic given this initial choice.

\paragraph{Stretchable sorting networks.} The following is one way to
generate a sorting network. Consider a set of $n$ points in
general position in $\R^2$, and label them $1,\ldots,n$ in order
of increasing $x$-coordinate. Now rotate the set of points by an
angle $\theta$. For all but finitely many $\theta$, listing the
labels of the points in order of increasing $x$-coordinate gives a
permutation in ${\cal S}_n$. And if we increase $\theta$
continuously from $0$ to $\pi$, these permutations yield the
sequence of configurations for a sorting network. Not all sorting
networks can be obtained in this way; in fact those which can are
exactly those whose wiring diagram may be drawn in the plane so
that all the trajectories are straight lines; such networks are
called {\em stretchable}---see \cite{goodman} for details. (The
smallest non-stretchable network, unique up to symmetries, is the
$5$-element example $\omega=(1,3,4,2,1,3,4,2,1,3)$.) In the proof
of \thmref{great-circle} we will see that the assumption of that
theorem implies that the sorting network is approximated by a
stretchable network obtained by rotating a set of points in $\R^2$
which approximate the Archimedes measure $\arch_{1/2}$.

Consider an $n$-element USN, and choose $m$ out of the $n$
particles uniformly at random, independently of the USN.  If we
observe only the relative order of these $m$ particles then we
obtain a random $m$-element sorting network.  If Conjecture
\ref{conj-great-circle} holds then it may be deduced that, as
$n\to\infty$ with $m$ fixed, the distribution of this sorting
network converges to a measure whose support is exactly the set of
stretchable $m$-element networks.  This follows from the proofs in
Section \ref{S:disc}.

\paragraph{Gallery.} For more simulation pictures, see the gallery\\
\verb+http://www.math.ubc.ca/~holroyd/sort+.

\section{Preliminaries}

In this section we present some definitions and basic results.
\begin{proof}[Proof of \thmref{stat-semi}(i)]
If $\omega=(s_1,\ldots ,s_N)$ is any sorting network then it is easily seen that
$$\omega':=(s_2,\ldots, s_N, n-s_1)$$
is also a sorting network, and furthermore that the map
$\omega\mapsto \omega'$ is a bijection from $\sort_n$ to
$\sort_n$.  The result follows immediately.
\end{proof}

We note also that
\begin{equation} \label{sym1}
(s_1,\ldots ,s_N)\mapsto(s_N,\ldots ,s_1)
\end{equation}
 and
\begin{equation}\label{sym2}
(s_1,\ldots ,s_N)\mapsto(n-s_1,\ldots ,n-s_N)
\end{equation}
 are bijections from
$\sort_n$ to $\sort_n$, so the measure $\pusn^n$ has the
corresponding symmetries.

For a permutation $\sigma\in\ess_n$, denote the {\bf inversion
number}
\[
  \inv(\sigma) = \# \Big\{(i,j):\;1 \le i < j \le n\text{ and }
  \sigma(i)>\sigma(j)\Big\}.
\]
It is straightforward to see that $\inv(\sigma)$ is the
graph-theoretic distance from the identity to $\sigma$ in the
Cayley graph of ${\cal S}_n$ generated by the swaps
$\{\tau_1,\ldots,$ $\tau_{n-1}\}$. Hence in any sorting network we
have $\inv(\sigma_k)=k$ for all $k$.

\section{Young tableaux} \label{Young tableaux}

Young tableaux are a central tool in our proofs; we start by
introducing some standard notation and facts. Let $N\in
\mathbb{N}:=\{1,2,\ldots\}$. A {\bf partition of $N$} is a
sequence $\lambda = (\lambda_1, \lambda_2, \ldots, \lambda_k)$ of
positive integers such that $\lambda_1 \ge \lambda_2 \ge \ldots
\ge \lambda_k$ and $N =\sum_i \lambda_i$. We denote
$|\lambda|:=N$. We identify each partition $\lambda$ with its
associated {\bf Young diagram}, which is the set $\{(i,j)\in
\mathbb{N}^2 \ :\ 1\le i\le k,\ 1\le j\le \lambda_i\}$.
Traditionally each element $(i,j)$ (called a {\bf cell}) in the
diagram is drawn as a square, in the coordinate system with
$(1,1)$ at the top-left and $(1,2)$ to its right.  We denote the
set of partitions of $N$ by $\partit(N)$.

Two Young diagrams will play a central role: the {\bf $n\times n$
square diagram} $(n,n,\ldots,n)$, which we denote by $\square_n$,
and the {\bf staircase diagram} $(n-1,n-2,\ldots,1)$, which we
denote by $\staircase_n$.

If $\lambda\in\partit(N)$, let $\lambda'=(\lambda_1',\lambda_2',
\ldots,\lambda_d')$ denote the {\bf conjugate partition} to $\lambda$,
where $d=\lambda_1$ and $\lambda_i' = \#\{ 1\le j\le k\ :\ \lambda_j
\ge i \}$. The conjugate partition corresponds to the Young diagram
obtained by reflecting the Young diagram of $\lambda$ along the
northwest-southeast diagonal.

\begin{figure}
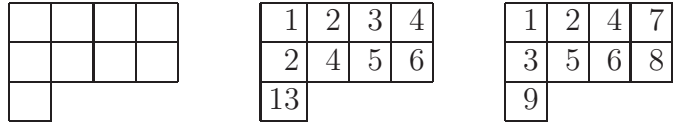

\begin{center}
  \begin{tabular}{|r|r|r|r|}
    \hline
    \ \  &\ \ &\ \ &\ \ \\ \hline
    & & & \\ \hline
    \\ \cline{1-1}
  \end{tabular}
  \qquad
  \setlength{\tabcolsep}{3pt}
  \begin{tabular}{|r|r|r|r|}
    \hline
    \ 1 & \ 2 &\ 3 &\  4 \\ \hline
    2 & 4 & 5 & 6 \\ \hline
    13 \\ \cline{1-1}
  \end{tabular}
  \qquad
  \begin{tabular}{|r|r|r|r|}
    \hline
    \ 1 &\ 2 &\ 4 &\ 7 \\ \hline
    3 & 5 & 6 & 8 \\ \hline
    9 \\ \cline{1-1}
  \end{tabular}
\end{center}
\caption{The Young diagram $(4,4,1)$, a Young tableau and a standard Young
  tableau.}
\label{F:SYT}
\end{figure}

A {\bf Young tableau} of shape $\lambda$, where $\lambda \in
\partit(N)$, is an assignment of positive integers, called
{\bf entries}, to the cells of $\lambda$ such that every row and
column of the diagram contain increasing sequences of numbers. A
{\bf standard Young tableau (SYT)} is a Young tableau in which the
numbers assigned to all the cells are $1,2,\dots,N$.  See
\figref{SYT}. We denote the set of SYT of shape $\lambda$ by
$\syt(\lambda)$, and we denote $d(\lambda) = \#\syt(\lambda)$
(sometimes called the \emph{dimension} of $\lambda$ in
representation-theoretic contexts), the number of standard Young
tableaux of shape $\lambda$.  Frame, Robinson and Thrall
\cite{framerobinsonthrall}, \cite[Sec.\ 1.5.4]{knuth} proved the
following formula for $d(\lambda)$.

\begin{thm}[Hook formula; Frame, Robinson and Thrall]
\label{hook}
For each cell $(i,j)\in\lambda$ let
$h_{i,j}(\lambda):=\lambda_i-j+\lambda_j'-i+1$ be the {\bf hook
number} of $(i,j)$ in $\lambda$. Then
$$d(\lambda) = \frac{|\lambda| !}{\prod_{(i,j)\in\lambda} h_{i,j}(\lambda)}. $$
\end{thm}

For two Young diagrams $\mu,\lambda$ write $\mu \nearrow \lambda$
(``$\mu$ increases to $\lambda$'') to mean that $\lambda$ can be
obtained from $\mu$ by the addition of one cell. The {\bf Young
lattice} is the directed graph whose vertex set is
$\cup_{N=0}^\infty \text{Par}(N)$
and whose edges are all the pairs $(\mu,\lambda)$ with
$\mu\nearrow \lambda$. Standard Young tableaux of shape $\lambda$
are in bijection with paths in the Young lattice leading from the
empty diagram $\emptyset$ to $\lambda$: to the path $\emptyset =
\lambda_0\nearrow \lambda_1\nearrow \lambda_2\nearrow \cdots
\nearrow \lambda_N=\lambda$ we attach the SYT which records the
order in which new cells were added to the diagrams along the
path, i.e., the unique tableau $T=(t_{i,j})_{(i,j)\in\lambda}$
such that for all $0\le k\le N$ we have that
$$ \lambda_k = \{ (i,j)\in \lambda : t_{i,j} \le k \}. $$
We call $T$ the {\bf recording tableau} of the increasing sequence of
diagrams $(\lambda_k)_{0\le k\le N}$.

As an illustration of the use of the hook formula for sorting
networks, we verify the claim of Example \ref{double-flip}. We
will use the fact (see \cite[p.\ 135]{finch}) that
\[
  J(n) := \prod_{j=1}^n j^j = \exp\left[\frac{n^2+n}{2}\log n +
  O(n)\right] \ntoinf.
\]
\nopagebreak Thus, using the hook formula and Stirling's formula
we can compute
\begin{eqnarray}
\nonumber
 d(\staircase_n)
 &=&
 \frac
   {\binom{n}{2}!}
   {1^{n-1}3^{n-2}5^{n-3}\cdots (2n-3)^1}
 =
 \frac
   {\binom{n}{2}!\left(\frac{J(2n-2)}{2^{n(n-1)}
       J(n-1)^2}\right)^{1/2}}
   {\left( \frac{(2n-2)!}{2^{n-1}(n-1)!}\right)^{n-\frac{1}{2}}}
 \\
 \label{eq:num-tri}
 &=&
 \exp \left[\frac{n^2-n}{2}\log n +
   (\tfrac{1}{4}-\log 2)n^2 + O(n)\right];
\end{eqnarray}
\begin{eqnarray}
 \nonumber
 d(\square_n)
 &=&
  \frac
   {\left(n^2\right)!}
   {J(n) (n+1)^{n-1}(n+2)^{n-2}\cdots (2n-1)^1}
 \\& =&
 \label{eq:num-sq}
 \frac
   {\left(n^2\right)! \frac{J(2n-1)}{J(n)}}
   {J(n)\left[\frac{(2n-1)!}{n!}\right]^{2n}}
 = \exp \Big[n^2\log n
  + (\tfrac{1}{2}-\log 4)n^2 + O(n)\Big].
\end{eqnarray}
By \eqref{eq:num-tri} and \eqref{eq:numsortnet} we have $\#
\sort_n=d(\staircase_n)$ (also see Section \ref{sec:eg} below).

For a permutation $\nu \in S_n$, a {\bf partial sorting network}
(also called a reduced word) of $\nu$ is a sequence
$(s_1,s_2,\ldots,s_k)$ such that $\nu = \tau_{s_1}\tau_{s_2}\cdots
\tau_{s_k}$ and $k=\inv(\nu)$. Let $R(\nu)$ denote the number of
partial sorting networks of $\nu$. In general, evaluation of
$R(\nu)$ is a deep problem---see e.g.\ \cite{garsia,little}.

Let $\nu$ be any permutation, and let $k=\inv(\nu)$. Then
$\omega=(s_1,s_2,\ldots,s_N)$ is a sorting network passing through
configuration $\nu$ if and only if $(s_1,\ldots,s_k)$ is a partial
sorting network for $\nu$ and $(s_{k+1},\ldots,s_N)$ is a partial
sorting network for $\nu^{-1}\rho$.  Hence the probability that
the USN passes through $\nu$ equals
\begin{equation}
\label{pass-through}
\pusn(\sigma_k=\nu)=\frac{R(\nu)R(\nu^{-1}\rho)}{R(\rho)}.
\end{equation}

\begin{proof}[Proof of Example \ref{double-flip}]
  In the case of $\psi$, we can compute the factors in
  \eqref{pass-through} above explicitly. We have
  $\inv(\psi)=h$. Firstly, $R(\psi)$ is equal to $\binom{h}{h/2}
  d(\staircase_{n/2})^2$, since to get from
  $\id$ to $\psi$ one must reverse the particles
  $1,\dots,\frac{n}{2}$, and independently reverse the particles
  $\frac{n}{2}+1,\dots,n$, with $\binom{h}{h/2}$ choices for the order
  in which to intersperse the left- and the right-half swaps.

  Secondly, we claim that the number of partial sorting networks of
  $\psi^{-1}\rho = (\frac{n}{2}+1,\dots,n,\; 1,\dots,\frac{n}{2})$ is
  equal to $d(\square_{n/2})$. This is because, given such a partial
  sorting network $(s_1,s_2,\ldots,s_{N-h})$, we can construct a
  standard Young tableau of shape $\square_{n/2}$ whose $i$th row
  lists the times $k_1 < k_2 < \ldots < k_{n/2}$ at which particle $i$
  moved, and it is easy to see that this map is a bijection from the
  set of partial sorting networks of $\psi^{-1}\rho$ onto
  $\syt(\square_{n/2})$.
Thus we have:
\begin{eqnarray*}
     R(\psi)&=&{\textstyle \binom{h}{h/2}} d(\staircase_{n/2})^2
            = 2^{n^2/4 + O(n)}  d(\staircase_{n/2})^2; \\
     R(\psi^{-1}\rho)&=&d(\square_{n/2});  \\
     R(\rho)&=&d(\staircase_{n}).
\end{eqnarray*}
An application of the asymptotics (\ref{eq:num-tri}) and
(\ref{eq:num-sq})  for the number of tableaux together with
(\ref{pass-through}) verifies the claim of Example
\ref{double-flip}.  Interestingly, the leading terms in
$n^2\log n$ cancel in the exponent.
\end{proof}

\section{The Edelman-Greene bijection}
\label{sec:eg}

Stanley, who proved \eqref{eq:numsortnet}, noticed that by the
hook formula the right-hand side of \eqref{eq:numsortnet} is equal
to $d(\staircase_n)$, the number of staircase shape standard Young
tableaux of order $n$. Later, Edelman and Greene
\cite{edelmangreene} found an explicit bijection between
$\syt(\staircase_n)$ and $\sort_n$. This bijection will play an
important part in what follows, so we describe it and its inverse
now.

Given a standard Young tableau $T\in\syt(\lambda)$,
where $N=|\lambda|$,
denote by $(i_{\max}(T),j_{\max}(T))$ the coordinates of the cell containing
the maximum entry $N$ in $T$.

\begin{figure}
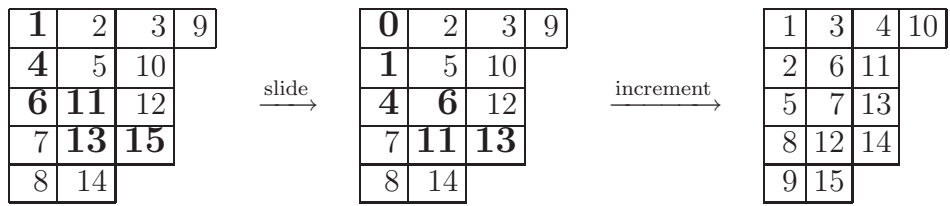

\begin{center}
  \setlength{\tabcolsep}{3pt}
  \begin{tabular}{|r|r|r|r|}
    \hline
    \ \bf{\large 1} &      2  &      3  & \ 9 \\ \hline
      \bf{\large4} &      5  &     10  \\ \cline{1-3}
      \bf{\large6} & \bf{\large11} &     12  \\ \cline{1-3}
          7  & \bf{\large13} & \bf{\large15} \\ \cline{1-3}
          8  &     14 \\ \cline{1-2}
  \end{tabular}
  \quad $\xrightarrow[]{\text{slide}}$ \quad
  \begin{tabular}{|r|r|r|r|}
    \hline
    \ \bf{\large0} &      2  &      3 & \ 9 \\ \hline
      \bf{\large1} &      5  &     10 \\ \cline{1-3}
      \bf{\large4} &  \bf{\large6} &     12 \\ \cline{1-3}
          7  & \bf{\large11} & \bf{\large13} \\ \cline{1-3}
          8  &     14 \\ \cline{1-2}
  \end{tabular}
  \quad $\xrightarrow[]{\text{increment}}$ \quad
  \begin{tabular}{|r|r|r|r|}
    \hline
    \ 1 &  3 &  4 & 10 \\ \hline
      2 &  6 & 11 \\ \cline{1-3}
      5 &  7 & 13 \\ \cline{1-3}
      8 & 12 & 14 \\ \cline{1-3}
      9 & 15 \\ \cline{1-2}
  \end{tabular}
\end{center}
\caption{The sliding sequence and the Sch\"utzenberger operator.
Shown are: (a) A tableau $T$. In bold is the sliding sequence
(obtained by starting from the maximum entry and repeatedly
passing to the larger of the entries above and left); (b) the
tableau obtained by sliding the entries down along the sliding
sequence; (c) the tableau $\evac(T)$.} \label{F:evac}
\end{figure}
Define the {\bf Sch\"utzenberger operator} $\evac:\syt(\lambda)
\to \syt(\lambda)$ as follows.  Start with a tableau
$T=(t_{i,j})_{(i,j)\in \lambda}$. Construct the {\bf sliding
sequence} of cells $c_0, c_1,\ldots,c_d \in \lambda$, where $c_0 =
(i_{\max}(T),j_{\max}(T))$ and $c_d=(1,1)$, by the requirements
that $c_r-c_{r+1} = (1,0)\textrm{ or }(0,1)$ for all $0\le r\le
d-1$, and $c_r-c_{r+1} = (1,0)$ if and only if
$t_{c_r-(1,0)}>t_{c_r-(0,1)}$ (where we adopt the notational
convention that for a cell $(i,j)$ with either of $i,j$ being
non-positive we have $t_{i,j}=-\infty$). Then the tableau
$\evac(T)=(t_{i,j}')_{(i,j)}$ is defined by setting
$t_{c_r}'=t_{c_{r+1}}+1$ for $0\le r\le d-1$, $t_{1,1}'$=1, and
$t_{i,j}'=t_{i,j}+1$ for all other cells $(i,j)\in\lambda$. The
definition is illustrated in \figref{evac}. It is easy to see that
$\evac$ is a bijection of $\syt(\lambda)$ onto itself.

\begin{trivlist}
\item
{\bf Definition.} \em
  The Edelman-Greene bijection $\EG:\syt(\staircase_n) \to\sort_n$ is
  defined by
  \[
  \EG(T) = \Big(j_{\max}(\evac^{N-k}(T))\Big)_{k=1,\dots,N}
  \]
  where as before $N=\binom{n}{2}$, and $\evac^k$ denotes the $k$th
  iterate of $\evac$.
\end{trivlist}

It is far from obvious that the map $\EG$ is a bijection to
$\Omega_n$, nor what its inverse looks like. It turns out that the
inverse may be described in terms of a Young tableau construction
algorithm which is a modification of the RSK algorithm (see
\cite[Ch.\ 7.11]{stanley}). Given a sorting network
$\omega=(s_1,s_2,\ldots,s_N) \in\sort_n$, we construct a sequence
of (non-standard) Young tableaux $T_0, T_1, \ldots, T_N$ whose
shapes $\emptyset = \lambda_0, \lambda_1, \ldots,
\lambda_N=\staircase_N$ form an increasing sequence of diagrams,
i.e., $\lambda_i \nearrow\lambda_{i+1}$. To get $T_{i+1}$ from
$T_i$, apply the following {\bf insertion algorithm} to the input
$(T_i, s_{i+1})$.

\begin{figure}
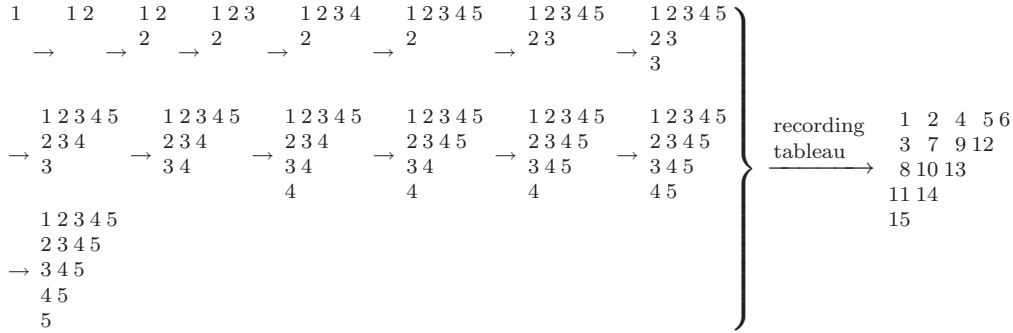

\begin{center}
  \def\YT#1{\begin{tabular}{rrrrr} #1 \end{tabular}}
  \setlength{\tabcolsep}{1pt}
  $\left.
  \begin{minipage}{.7\textwidth}
    \scriptsize
    \YT{1         \\         \\       \\     \\} $\to$
    \YT{1&2       \\         \\       \\     \\} $\to$
    \YT{1&2       \\ 2       \\       \\     \\} $\to$
    \YT{1&2&3     \\ 2       \\       \\     \\} $\to$
    \YT{1&2&3&4   \\ 2       \\       \\     \\} $\to$
    \YT{1&2&3&4&5 \\ 2       \\       \\     \\} $\to$
    \YT{1&2&3&4&5 \\ 2&3     \\       \\     \\} $\to$
    \YT{1&2&3&4&5 \\ 2&3     \\ 3     \\     \\} $\to$
    \YT{1&2&3&4&5 \\ 2&3&4   \\ 3     \\     \\} $\to$
    \YT{1&2&3&4&5 \\ 2&3&4   \\ 3&4   \\     \\} $\to$
    \YT{1&2&3&4&5 \\ 2&3&4   \\ 3&4   \\ 4   \\} $\to$
    \YT{1&2&3&4&5 \\ 2&3&4&5 \\ 3&4   \\ 4   \\} $\to$
    \YT{1&2&3&4&5 \\ 2&3&4&5 \\ 3&4&5 \\ 4   \\} $\to$
    \YT{1&2&3&4&5 \\ 2&3&4&5 \\ 3&4&5 \\ 4&5 \\} $\to$
    \YT{1&2&3&4&5 \\ 2&3&4&5 \\ 3&4&5 \\ 4&5 \\ 5}
  \end{minipage}
  \right\} \xrightarrow[]{\text{\shortstack[l]{recording\\tableau}}}$
  {\scriptsize
    \YT{1&2&4&5&6 \\ 3&7&9&12 \\ 8&10&13 \\ 11&14 \\ 15}
  }
  \caption{Computation of $\EG^{-1}(1,2,1,3,4,5,2,1,3,2,1,4,3,2,1)$.}
  \label{F:eg-inv}
\end{center}
\end{figure}
\paragraph{Insertion algorithm.} {\it
Given a Young tableau $T=(t_{i,j})_{(i,j)\in\lambda}$ of shape $\lambda$
and a positive number $u$, construct a new tableau $T'=(t_{i,j}')$ whose
shape is the union of $\lambda$ with one new cell, as follows.}
\begin{alglist}
\item
{\bf (Initialize).}
 Set $k\leftarrow 1$ and $q \leftarrow
u$. Set $t_{i,j}' \leftarrow t_{i,j}$ for all
$(i,j)\in\mathbb{N}^2$, with the convention that $t_{i,j}=\infty$
for a cell $(i,j)\in \mathbb{N}^2\setminus\lambda$.
\item
{\bf (Find next bumping cell).} Set $\ell$ to be the least positive
integer $j$ such that $t_{k,j} \ge q$. Set $t_{k,\ell}' \leftarrow q$.
If $q=t_{k\ell}$, set $q\leftarrow q+1$, otherwise set $q\leftarrow t_{k\ell}$.
Set $k\leftarrow k+1$.
\item
If $q=\infty$, terminate and return the enlarged tableau $T'$.
Otherwise return to Step 2.
\end{alglist}

\begin{trivlist}
\item
{\bf Definition.} \em
  The inverse Edelman-Greene bijection $\EG^{-1}:\sort_n\to
  \syt(\staircase_n)$ is defined by setting $\EG^{-1}(\omega)$ to be the
  recording tableau of the sequence of Young diagrams $\lambda_0\nearrow
  \lambda_1\nearrow \dots \nearrow \lambda_N=\staircase_n$ constructed
  above.
\end{trivlist}

\figref{eg-inv} shows $\EG^{-1}$ applied to the sorting network of
\figref{wire}. The following theorem justifies these definitions. The proof
can be found in \cite{edelmangreene}; see also \cite{felsner}.

\begin{thm}[Edelman \hspace{-1mm} and \hspace{-1mm} Greene] \label{T:EG_bij}
  The map $\EG$ is a bijection from $\syt(\staircase_n)$ to $\sort_n$, and
  the map $\EG^{-1}$ is its inverse.
\end{thm}

As a first application of the Edelman-Greene bijection, we prove
an exact formula for the distribution of the first swap location
$s_1 = s_1(\omega_n)$ of a uniform $n$-element sorting network
$\omega_n$, and use it to prove \thmref{stat-semi}(ii).

\begin{prop}[Swap distribution] \label{semi}
If $\omega_n$ is a uniform $n$-element sorting network, then
\begin{equation}
\pusn^n( s_1 = r ) = \frac{1}{N}\cdot
\frac{\big(3\cdot 5\cdot 7\cdots (2r-1)\big)
  \big(3\cdot 5\cdots (2(n-r)-1)\big)}
{\big(2 \cdot 4\cdot 6\cdots (2r-2)\big)
 \big(2 \cdot 4\cdots (2(n-r)-2)\big)}.
\end{equation}
\end{prop}

\begin{proof}
Let $1\le r\le n-1$. By the definition of $\EG$, the sorting
networks $\omega =(s_1,s_2,\ldots,s_N)\in\sort_n$ for which $s_1 =
r$ are exactly the ones for which the standard Young tableau
$\EG^{-1}(\omega)$ has its maximum entry in the cell $(n-r,r)$. Since
$\EG$ is a bijection, the number of such $\omega$'s is the number
of SYTs of shape $\staircase_n\setminus\{(n-r,r)\}$. Thus
\[
\pusn^n(s_1 = r) =
\frac{d(\staircase_n\setminus\{(n-r,r)\})}{d(\staircase_n)}.
\]

Write this using the hook formula, Theorem \ref{hook}, to yield the result.
\end{proof}

\begin{proof}[Proof of \thmref{stat-semi}(ii)]
Denote $a_{n,r} = \pusn^n(s_1 =r), 1\le r\le n-1$. Observe that by
Proposition \ref{semi}  we have
\begin{eqnarray*} a_{n,r} &=& \frac{2}{n(n-1)} \cdot
\frac{(2r)(2r)!}{2^{2r} (r!)^2}\cdot
\frac{2(n-r)(2(n-r))!}{2^{2(n-r)} ((n-r)!)^2}
\\ &=& \frac{8\sqrt{r(n-r)}}
{\pi n(n-1)}\cdot \frac{\sqrt{\pi r}\binom{2r}{r}}{2^{2r}}
\cdot \frac{\sqrt{\pi (n-r)} \binom{2(n-r)}{n-r}}{2^{2(n-r)}}
\end{eqnarray*}
Therefore, using Stirling's formula in its explicit form $ 1 \le
m!(2\pi m)^{-1/2}\left(\frac{e}{m}\right)^m$ $\le
1+\frac{1}{12m-1} $ (an immediate consequence of \cite[Eq (9.15),
p.\ 54]{feller}), we get that
\begin{eqnarray*}
\left(1-\frac{1}{6r}\right)\left(1-\frac{1}{6(n-r)}\right) &\le&
\frac{a_{n,r}}{\frac{8}{\pi n(n-1)} \sqrt{r(n-r)}} \\ &\le&
\left(1+\frac{1}{24r-1}\right) \left(1+\frac{1}{24(n-r)-1}\right).
\end{eqnarray*}
This implies easily that for all $-1<a<b<1$ we have
\[
\pusn^n\left(a \le \frac{2s_1}{n} - 1 \le b\right)
= \sum_{\frac{n}{2}(a+1) \le r \le \frac{n}{2}(b+1)} a_{n,r}
\xrightarrow[n\to\infty]{} \frac{2}{\pi}\int_a^b \sqrt{1-t^2}dt,
\]
as required.
\end{proof}

\section{Limit profile for staircase tableaux}
\label{sec:limit}

For any Young diagram $\lambda$ we write $\P_\lambda$ for the
uniform measure on the set $SYT(\lambda)$ of standard Young
tableaux. It is natural to consider the limiting behaviour of a
random tableau of distribution $\P_{\lambda_n}$ for a sequence of
diagrams $(\lambda_n)$ of a given shape and increasing size.  For
general shape, the problem of rigorously determining the complete
limiting profile is open (see, however \cite{Ke1,Ke2} and
\cite[Theorem 1.5.1]{Biane}).  An exception is the square diagram
$\square_n$, where the problem was solved by Pittel and Romik
\cite{pittelromik}.  In this section we use their result to derive
a solution for the staircase diagram $\staircase_n$.

We start by stating the main result from \cite{pittelromik}.  It
will be convenient to use the following coordinate system. If
$(i,j)$ is a cell of $\square_n$, then its rotated (and scaled)
coordinates are
\begin{align*}
  u&= u(i,j) := \frac{i-j}{n};   &   v& = v(i,j) := \frac{i+j}{n},
\end{align*}
(note that this differs from the coordinate system in
\cite{pittelromik} by a factor of $\surd2$).

\begin{figure}
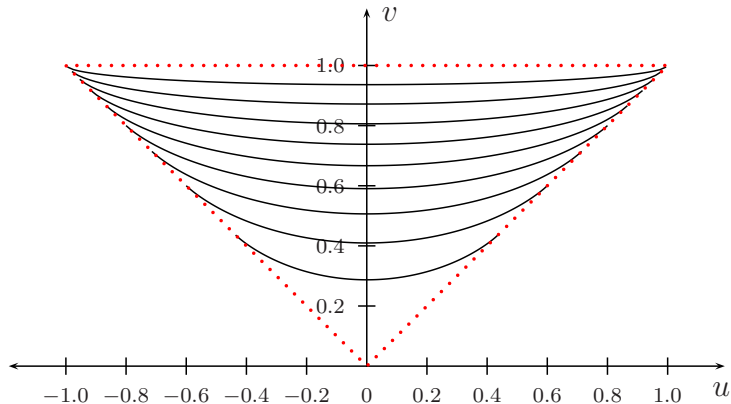

\begin{center}
    \psset{unit=4cm,linewidth=0}
    \pspicture(-1,-.1)(1,1)
    {\scriptsize
    \psaxes[linewidth=0.005,Dx=.2,Dy=.2]{<->}(0,0)(-1.19,0)(1.19,1.19) }
    \psset{plotpoints=100}
    \multido{\n=0.1+.1}{9}{
      \parametricplot[linewidth=0.005]{-1}{1}{
        /beta {\n} def
        /x 1 beta beta mul sub sqrt t mul def
        x
        /Q 1 beta beta mul sub x x mul sub def
        /R Q 0 gt {Q sqrt beta div} {.000001} ifelse def
        x R div ATAN x mul R ATAN add 2 PI div mul}
    }
    \psline[linecolor=red,linewidth=0.012,linestyle=dotted](0,0)(1,1)(-1,1)(0,0)
    \put(1.15,-.1){$u$}
    \put(0.05,1.15){$v$}
  \endpspicture
  \caption{The curves $v=h_\alpha(u)$ for $\alpha=0.1, 0.2, 0.3, \ldots,
    0.9$, bounded between the graphs of $v=|u|$ and $v=1$.}
  \label{F:levcurves}
\end{center}
\end{figure}

We define the following functions, which will describe the
limiting profile.  For $\alpha\in[0,2]$ the function $h_\alpha :
[-\sqrt{\alpha(2-\alpha)},\sqrt{\alpha(2-\alpha)}]\to[0,1]$ is
defined by
\begin{equation} \label{eq:hdef}
  h_\alpha(u) := \frac{2}{\pi}\left[
    u \arctan\left( \frac{u}{R} \right) + \arctan^{-1} R\right]
  \quad\text{where}\quad
  R = \frac{\sqrt{\alpha(2-\alpha)-u^2}}{1-\alpha},
\end{equation}
for $\alpha\in[0,1]$ (where $\tan^{-1}\infty:=\pi/2$ giving $h_1\equiv 1$),
and by
\[
h_{2-\alpha}(u):=2-h_\alpha(u)
\]
for $\alpha\in(1,2]$. The curve $v=h_\alpha(u)$ will approximate the
level-$(\alpha n^2 /2)$ contour of the tableau; \figref{levcurves} shows
some of these curves. The function $L:[0,1]\times[0,1]\to[0,2]$ is defined
implicitly by
\[
L\left(\frac{u+v}{2},\frac{v-u}{2}\right) = \alpha  \iff  h_\alpha(u) = v.
\]

The following result \cite[Theorem 1(i)]{pittelromik}) gives the
limiting profile for uniform square Young tableaux.

\begin{thm}[Limit profile for square tableaux; Pittel and Romik]
  \label{T:square_limit}
  Let $\P_{\square_n}$ be the uniform measure on Young tableaux
  $(s_{i,j})_{(i,j)\in\square_n} \in \syt(\square_n)$. For any
  $\eps>0$,
  \[
  \P_{\square_n}\left( \max_{(i,j)\in\square_n} \left| \frac{2 s_{i,j}}{n^2} -
                      L\left(\frac{i}{n},\frac{j}{n}\right) \right| > \eps
                  \right)
  \xrightarrow[n\to\infty]{} 0.
  \]
\end{thm}

We shall deduce the following analogous result for the limit profile of a
staircase tableau, where the function $L$ is the same as above.

\begin{thm}[Limit profile for staircase tableaux] \label{T:staircase_limit}
  Let $\P_{\staircase_n}$ be the uniform measure on Young tableaux
  $(t_{i,j})_{(i,j)\in\staircase_n} \in \syt(\staircase_n)$. For any
  $\eps>0$,
  \[
  \P_{\staircase_n}\left( \max_{(i,j)\in\staircase_n} \left| \frac
    {2 t_{i,j}}{n^2} - L\left(\frac{i}{n},\frac{j}{n}\right) \right| > \eps
                  \right)
  \xrightarrow[n\to\infty]{} 0.
  \]
\end{thm}

Thus the limit profile for the staircase tableau is the same as that for
half of the square tableau. Other results in \cite{pittelromik} give
explicit bounds on deviations from the limit profile, but only in the
interior of the square. These estimates may be translated to staircase
tableaux as well. However, uniform convergence in probability is sufficient
for our purposes. It is important that \thmref{staircase_limit} includes
the boundary of the diagram.

Our main tool in proving the above is the following general result
concerning continuity of random tableaux in the shape. For Young
diagrams $\lambda,\mu$ we write $\lambda\subseteq\mu$ if this
relation holds for $\lambda,\mu$ as subsets of ${\mathbb N}^2$.

\begin{thm}[Coupling] \label{T:SYT_couple}
  Let $\lambda\subseteq\mu$ be a pair of Young diagrams. There exists a
  coupling of the measures $\P_\mu$ on $S=(s_{i,j})\in\syt(\mu)$ and
  $\P_\lambda$ on $T=(t_{i,j})\in\syt(\lambda)$ such that for all
  $(i,j)\in\lambda$
  \[
    s_{i,j} \le t_{i,j} + |\mu\setminus\lambda|.
  \]
\end{thm}

To prove \thmref{SYT_couple} we will make use of an algorithm from
\cite{greeneetal} for sampling from $\P_\lambda$. First note that, in
order to simulate a tableau with distribution $\P_\lambda$, it
suffices to be able to choose the location $\cmax=(i_\text{\rm
max},j_\text{\rm max})$ of the maximum entry $|\lambda|$ with the
correct distribution.  For then, after inserting this entry, we
may iteratively apply the same algorithm to the smaller diagram
$\lambda\setminus \{c_{\text{\rm max}}\}$ to locate the second largest
entry, and so on.

For a cell $(i,j)\in\lambda$, define its {\bf hook} to be the set
\[
H_{(i,j)}(\lambda):=\big\{(k,j)\in\lambda: k\ge i\big\}\cup
\big\{(i,k)\in\lambda: k\ge j\big\}.
\]
The location $\cmax$ of the maximum entry may be simulated using
the following {\bf hook walk} algorithm from \cite{greeneetal}.

\paragraph{Hook walk algorithm.}  {\em Given a Young diagram $\lambda$,
choose a random sequence of cells $c_0,\ldots,c_r$ iteratively as follows.}
\begin{alglist}
\item Choose a cell $c_0$ uniformly at random from $\lambda$.
\item Given that cells $c_0,\ldots, c_{k-1}$ have been chosen, choose $c_k$
  uniformly at random from the hook $H_{c_{k-1}}(\lambda)$.
\item Repeat Step 2 until we obtain a cell $c_r$ with $\#
  H_{c_r}(\lambda)=1$, then stop.
\end{alglist}

\begin{thm}[Hook walk; Greene, Nijenhuis and Wilf]
\label{hookwalk}
  The random final cell $c_r$ constructed by the hook walk has the same
  distribution as $\cmax$ under $\P_\lambda$.
\end{thm}

\begin{lemma}[Domination] \label{L:one_step}
  Assume $\lambda\subseteq\mu$, and let $\lambda':=\lambda\setminus
  \{\cmax(\lambda)\}$ and $\mu':=\mu\setminus \{\cmax(\mu)\}$ be the
  random Young diagrams obtained by removing the largest entry in
  the respective uniform standard Young tableaux. Then
  we have the stochastic domination $\lambda' \subseteq_{\text{\rm st}} \mu'$.
\end{lemma}

\begin{proof}
  Consider the hook walk applied to $\lambda$ and $\mu$. It is enough to
  couple the two hook walks so that either they stop at the same cell, or
  the walk in $\mu$ stops at a cell in $\mu\setminus\lambda$. This will
  hold provided the two walks coincide until the first time the one in
  $\mu$ enters $\mu\setminus\lambda$. And this can be achieved as follows.
  Run the hook walk in $\mu$ according to the usual rules. Let the walk in
  $\lambda$ be identical to that in $\mu$ while the latter is in $\lambda$.
  If and when the walk in $\mu$ jumps to a cell in $\mu\setminus \lambda$,
  continue the walk in $\lambda$ according to the usual rules for $\lambda$
  using an independent source of randomness. It is easy to see that this
  gives the correct hook walk terminating probabilities for
  $c_{\max}(\lambda)$.  (A key observation is that a
  uniformly chosen element in the hook $H_{(i,j)}(\mu)$ conditioned to be
  in $H_{(i,j)}(\lambda)$ is distributed uniformly in $H_{(i,j)}(\lambda)$).
\end{proof}

\begin{proof}[Proof of \thmref{SYT_couple}]
  Construct the random tableaux $S,T$ iteratively by first choosing
  the maximum entry in each, then the second largest, and so on. Do
  this using the hook walks, and at each stage couple the two hook
  walks according to \lemref{one_step}, so that the remaining unfilled
  Young diagrams are always ordered. Let $m = |\mu\setminus\lambda|$.
  At the step when $k$ is entered into $\lambda$, say at location
  $(i,j)$, the entry $k+m$ is entered into $\mu$, and all subsequent
  entries entered into $\mu$ are $\le k+m$. By the ordering property
  one of those subsequent entries (possibly $k+m$)
  will be at the cell $(i,j)$.  Therefore $s_{i,j} \le k+m =
  t_{i,j}+m$.
\end{proof}

\begin{proof}[Proof of \thmref{staircase_limit}]
  Fix $\eps>0$, and consider a random square tableau $S=(s_{i,j})$ with law
  $\P_{\square_n}$. Let $\mu$ be the random Young diagram obtained by
  removing the cells with entries greater than $(1/2+\eps)n^2$ from
  $\square_n$,and define the event $A_n=\{\staircase_n\subseteq\mu\}$. From
  \thmref{square_limit} we find that $\P(A_n) \to 1$, because $L \equiv 1$
  along the diagonal.

  Note that, conditional on $\mu$, the tableau obtained by restricting $S$
  to $\mu$ has law $\P_\mu$. Also let $T=(t_{i,j})$ have law
  $\P_{\staircase_n}$. \thmref{SYT_couple} implies that $S$ and $T$ can be
  coupled so that on the event $A_n$ we have $s_{i,j}\le t_{i,j}+n^2\eps$
  for $(i,j)\in \staircase_n$. Thus by the square diagram result,
  \thmref{square_limit}:
\begin{eqnarray*}
  &&\P_{\staircase_n}\left( \max_{(i,j)\in\staircase_n}
   \left[-\frac{2 t_{i,j}}{n^2} +
                      L\left(\frac{i}{n},\frac{j}{n}\right) \right] > 3\eps
                  \right)
                  \le\\&&\hspace{3em}
  \P_{\square_n}\left( \max_{(i,j)\in\square_n}\left[
                    -\frac{2 s_{i,j}}{n^2} +
                      L\left(\frac{i}{n},\frac{j}{n}\right)\right]  > \eps
                  \right)+\P(A_n^c)
  \xrightarrow[n\to\infty]{} 0,
\end{eqnarray*}
as required. The corresponding upper bound is similar: instead we
begin by removing the entries greater than $(1/2-\eps)n^2$.
\end{proof}

We extract the following consequences of \thmref{staircase_limit}
for use in the later proofs.

\begin{cor}[First row] \label{C:first_row}
  For a staircase tableau $T=(t_{i,j})$ let $R_k = R_k(T) :=
  \max\{j:t_{1,j}\leq k\}$ be the number entries $\le k$ in the first row.
  For any $\eps>0$ we have
  \[
  \P_{\staircase_n}\big(\max_k |R_k-d_k|>\eps n\big)\xlim 0
  \]
  where
  $
  d_k := n \sqrt{\tfrac{k}{N}\left(2-\tfrac{k}{N}\right)}.
  $
\end{cor}

\begin{proof}
\thmref{staircase_limit} easily implies that for any $\delta$,
with probability tending to 1 as $n\to\infty$,
\begin{equation} \label{row_limit}
\max_j \left| \frac
    {t_{1,j}}{N} - L\left(0,\frac{j}{n}\right) \right| < \delta.
\end{equation}
From the definition of $L$, the map $x\mapsto L(0,x)$ is
continuous and strictly monotone, and satisfies
$L(0,\sqrt{\alpha(2-\alpha)})=\alpha$ for $\alpha\in[0,1]$.
Therefore given any $\epsilon>0$, we can choose $\delta>0$ such
that for all $x,\alpha\in[0,1]$,
$$|L(0,x)-\alpha|<\delta \text{ implies }
\big|x-\sqrt{\alpha(2-\alpha)}\big|<\eps.$$ We deduce that on the
event (\ref{row_limit}) we have
$$\max_j |j-d_{t_{1,j}}|<\eps n.$$
Since $k\mapsto d_k$ is strictly monotone this implies that
$\max_k |R_k-d_k|<\eps n$.
\end{proof}

\begin{cor}[Contours] \label{C:diag}
  Fix some $\alpha\in[0,1]$, and let ${\cal H}_\alpha={\cal H}_\alpha(T)$
  be the set of entries
  in a staircase tableau $T$ in the cells where $v>h_\alpha(u)$. For any
  $\eps>0$ we have the following bound on the symmetric difference:
  \[
  \P_{\staircase_n}\Big[\# \big({\cal H}_\alpha \diff \{\lceil \alpha N\rceil
  ,\dots,N\}\big) >
    \eps N \Big]
  \xrightarrow[n\to\infty]{} 0.
  \]
\end{cor}

\begin{proof}
Fix $\eps>0$.  By the continuity and strict monotonicity
of the function $L$, we may choose
$\delta>0$ such that the area of the region
$D:=\{(u,v):h_{\alpha-2\delta}(u)\le v\le h_{\alpha+2\delta}(u)\}$
is at most $\eps$.  Then for $n$ sufficiently large, on the event
$$\max_{(i,j)\in\staircase_n} \left| \frac
    {2 t_{i,j}}{n^2} - L\left(\frac{i}{n},\frac{j}{n}\right) \right|
    \leq\delta,$$
we have that all entries in the symmetric difference ${\cal
H}_\alpha \diff \{\lceil \alpha N\rceil ,\dots,N\}$ lie in $D$, so
the result follows from \thmref{staircase_limit}.
\end{proof}

\section{Law of large numbers}
\label{sec:LLN}

This section contains the proof of \thmref{LLN}. Recall the
semicircle measure $\semi(dx) = \frac2\pi \sqrt{1-x^2}\;
\ind_{x\in(-1,1)} \,dx$. Fix some interval $[a,b]\subset (-1,1)$,
and define for $0\leq s\leq t\leq 1$ and a sorting network
$\omega$:
\[
  S_{s,t}(\omega) := \# \Big\{ sN\le k<tN \ : \
                              \tfrac2n s_k(\omega)-1 \in [a,b] \Big\}.
\]

We will deduce \thmref{LLN} from the following.
\begin{lemma} \label{L:smalltime_lln}
  Fix an interval $[a,b]\subset (-1,1)$. For any $\delta$ small enough
  (depending on $a,b$),
  \[
  \pusn^n \Big( \big| S_{0,\delta} - \delta N \semi[a,b] \big| >
                8N\delta^2 \Big)
  \xrightarrow[n\to\infty]{} 0.
  \]
\end{lemma}

\begin{proof}[Proof of \thmref{LLN}]
  It suffices to prove that for any $[a, b]\subseteq[-1,1]$ and for
  any $\eps>0$ and $0\le s<t\le 1$ we have
  \begin{equation} \label{eq:lln_goal}
    \pusn^n \Big[ \Big| \tfrac1N S_{s,t} - (t-s)\semi[a,b] \Big| >
    \eps \Big] \xrightarrow[n\to\infty]{} 0.
  \end{equation}
  Since the total number of swaps is deterministically $N$, it is enough to
  prove this in the case $[a,b]\subset(-1,1)$. We deduce this from
  \lemref{smalltime_lln} as follows. Fix some positive integer
  $m$ to be chosen later,
  split the time interval $[s,t)$ into $m$ smaller intervals of length
  $\delta := \frac{t-s}{m}$, and define the events
  \[
  B_k := \left\{ \big| S_{s+k\delta,s+(k+1)\delta}(\omega)
                       - \delta N \semi[a,b] \big| > 9N\delta^2 \right\}
  , \qquad 0\le k<m.
  \]
  Let $B = \bigcup B_k$. By stationarity of the swap location process
  (\thmref{stat-semi}(i)), each of the random variables
  $S_{s+k\delta,s+(k+1)\delta}$ is within $\pm 1$ of a random variable having
  the same law as $S_{0,\delta}$.
  Hence by \lemref{smalltime_lln},
  \[
  \pusn^n(B) \le \sum_{k=0}^{m-1} \pusn^n(B_k) \xrightarrow[n\to\infty]{} 0.
  \]

  If $B$ does not occur then the quantity
  \[
  S_{s,t} = \sum_{k=0}^{m-1} S_{s+k\delta,s+(k+1)\delta}
  \]
  satisfies
  \[
  \big| S_{s,t} - m\delta N \semi[a,b] \big|  \le  9mN\delta^2;
  \]
  i.e., since $m\delta=(t-s)$
  \[
  \big| \tfrac1N S_{s,t} - (t-s) \semi[a,b] \big|
  \le 9(t-s)\delta \le 9\delta.
  \]
  Now \eqref{eq:lln_goal} follows by setting $m$ large enough that
  $9\delta<\eps$.
\end{proof}

\begin{proof}[Proof of \lemref{smalltime_lln}]
  By stationarity of the swap process, the first $\delta N$ swaps and the
  last $\delta N$ swaps have the same law. The idea of the proof
  is now as follows. From the
  Edelman-Greene bijection we see that the last $\delta N$ swaps are
  determined by the locations of the $\delta N$ largest entries in the
  staircase shaped Young tableau $T$ corresponding to $\omega$. By
  Corollary~\ref{C:diag}, the set of these locations is almost deterministic,
  which will imply our claim.

  For a Young tableau $T$, consider $j_{\max}(\evac^k T)$. To find it we
  start with the element $N-k$ in $T$, and perform $k$ iterations of
  $\evac$. At each iteration, the entry increases by 1, and possibly moves
  one square towards the diagonal. If it started close to the diagonal, it
  can only hit the diagonal in a limited region. In particular, if $N-k$
  started in region $A$ of \figref{LLN} then necessarily $\left(\frac2n
    j_{\max}(\evac^k T) - 1\right) \in [a,b]$. Similarly, if it started in
  either of the regions labelled $C$ then it will not exit through that
  interval. If $N-k$ started in the region marked $B$, then whether or not
  it exits in the interval $[a,b]$ depends on locations of other entries in
  the tableau.

  Let $\omega = EG(T)$, and let $A_\delta(T)$ be the number of entries
  greater than $(1-\delta)N$ in region $A$ of $T$, and similarly define
  $B_\delta(T)$ with region $B$. We find
  \begin{equation} \label{eq:SAB}
  0 \le S_{0,\delta}(\omega) - A_\delta(T) \le B_\delta(T).
  \end{equation}

  \begin{figure}
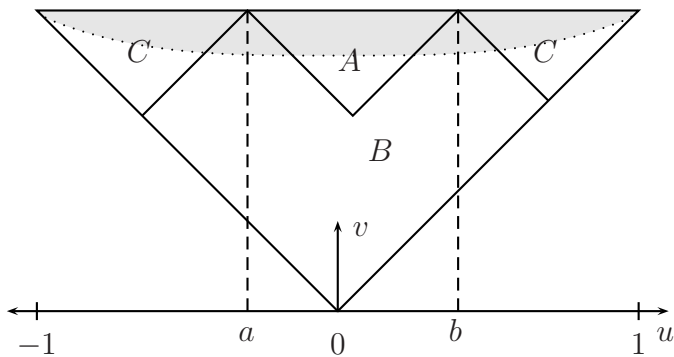

  \begin{center}
    \psset{unit=4cm}
    \pspicture(-1,-.1)(1,1)
    \psaxes[labels=x]{<->}(0,0)(-1.1,0)(1.1,.3)
    \newgray{verylight}{.9}
    \pscurve[linestyle=dotted,fillstyle=solid,fillcolor=verylight]
            (-1,1)(-.96,.97)(0,.85)(.96,.97)(1,1)
    \psline(0,0)(1,1)(-1,1)(0,0)
    \psline(-.65,.65)(-.3,1)(.05,.65)(.4,1)(.7,.7)
    \psline[linestyle=dashed](-.3,1)(-.3,0)
    \psline[linestyle=dashed](.4,0)(.4,1)
    \put(1.06,-.1){$u$}
    \put(0.05,.25){$v$}
    \put(-.33,-.1){$a$}
    \put(.37,-.1){$b$}
    \put(0,.8){$A$}
    \put(.1,.5){$B$}
    \put(-.7,.83){$C$}
    \put(.65,.83){$C$}
    \endpspicture
  \end{center}
  \caption{A Young tableau in the $(u,v)$ co-ordinate system. Entries in
    region $A$ can only exit through the interval $[a,b]$, while entries in
    $C$ cannot, and entries in $B$ may or may not. The region where
    $v>h_{1-\delta}(u)$ is shaded---this is the typical location of the
    entries greater than $(1-\delta)N$.}
  \label{F:LLN}
  \end{figure}

  To prove the lemma we show that with probability tending to 1, for a
  uniformly random tableau $T$, $B_\delta(T)$ and $\big| A_\delta(T) -
  \delta N \semi[a,b] \big|$ are both of order $\delta^2 N$. Here we use
  Corollary~\ref{C:diag}, with $\eps$ of the corollary equal to $\delta^2$.
  Consider the tableau $T$ in the $(u,v)$ co-ordinate system of Section
  \ref{sec:limit} as shown in \figref{LLN}. Let $H$ be the region where
  $v>h_{1-\delta}(u)$ for the function defined in \eqref{eq:hdef} (the
  shaded region in \figref{LLN}), and let ${\cal H}$ be the set of entries
  of $T$ in $H$.

  Corollary~\ref{C:diag} states that with probability tending to 1,
  \[
    \#\Big( {\cal H} \diff \{\lceil N-\delta N\rceil,\dots,N\} \Big)
    \le \delta^2 N.
  \]
  This implies that with probability tending to 1 we have
  \begin{align*}
    \Big|A_\delta(T) - N \lebb(H\cap A)\Big| &\le \delta^2 N, &
    \Big|B_\delta(T) - N \lebb(H\cap B)\Big| &\le \delta^2 N,
  \end{align*}
  where $\lebb$ denotes area of sets in the $(u,v)$
  plane.
  Applying these to \eqref{eq:SAB} we find that (with probability tending
  to 1)
  \begin{equation} \label{eq:Sleb}
  \Big| S_{0,\delta}(\omega) - N\lebb(H\cap A) \Big|
  \le \left( 2\delta^2 + \lebb(H\cap B) \right) N
  \end{equation}
  So, we need to estimate the areas of $H\cap A$ and $H\cap B$.

  We use the Taylor expansion of $h_{1-\delta}(u)$ around $\delta=0$, which is
  \[
    h_{1-\delta}(u) = 1 - \tfrac2\pi \sqrt{1-u^2} \, \delta +
    O(\delta^3)\quad\text{as }\delta\to 0,\text{uniformly in }u\in[a,b]
  \]
(see \cite[formula (7), p.\ 13]{pittelromik}).
  To estimate $\lebb(H\cap B)$ note that the side length of each of the two
  ``triangles'' comprising $H\cap B$ is of order $\delta$. Indeed, each is
  contained in a square with diagonal $\frac4\pi\delta$, and so
  \begin{equation} \label{eq:HB_small}
  \lebb(H\cap B) \le \frac{16}{\pi^2}\delta^2 < 2\delta^2.
  \end{equation}
  It remains to estimate $\lebb(H\cap A)$. The fact that $\left.
  \frac{\partial^2}{\partial\delta^2} h_{1-\delta}(u) \right|_{\delta=0}
  \equiv 0$ implies that for $\delta$ small enough (depending on $a,b$),
  for any $u\in[a,b]$ the error term in the Taylor expansion is at most
  $\delta^2$. Consequently, the area of $H\cap A$ can be estimated by
  integrating:
  \begin{equation} \label{eq:HA_int}
  \left| \lebb(H\cap A) - \delta \int_a^b \tfrac2\pi \sqrt{1-u^2} du \right|
  \le (b-a)\delta^2 + \lebb(H\cap B) \le 4\delta^2.
  \end{equation}
  (The $\lebb(H\cap B)$ term comes from the truncation near $a$ and $b$.)

  The result follows by applying \eqref{eq:HB_small} and \eqref{eq:HA_int} to
  \eqref{eq:Sleb}.
\end{proof}

\section{Octagon and H\"older bounds}

In this section we prove Theorems \ref{T:Holder} and
\ref{T:octagon}.
\begin{lemma}\label{det first row}
Let ${\mathcal R_k}={\mathcal R_k}(\omega)=(\lambda_k)_1$ be the
length of the first row of the Young diagram $\lambda_k$ created
by the first $k$ steps of the EG$^{-1}$ algorithm from a sorting
network $\omega$. Then $\sigma_k^{-1}(i)-i\le \mathcal R_k$ for
all $i,k$.
\end{lemma}

\begin{proof}

The first row of the recording and insertion tableaux during the
EG$^{-1}$ algorithm behave exactly the same way as during the
celebrated RSK algorithm. It is well known (see \cite{stanley})
that the RSK algorithm applied to any sequence of numbers creates
a tableau whose first row has length given by the longest
increasing subsequence (check this on Figure \ref{F:eg-inv}).

Here we only need the upper bound, which for completeness
we verify here.  Observe that the entries of the first row cannot increase as
the EG$^{-1}$ algorithm proceeds, as they only change through
bumping which replaces an element by something less or equal. So
given any increasing subsequence $a_1,\ldots, a_\ell$, we know
that each $a_i$ has to be inserted to the right of where $a_{i-1}$
was inserted. This shows that at any step $k$ the length of the
first row is at least the length of the longest increasing
subsequence of the input so far.

Now by time $k$ the position of particle $i$ has changed by
$\eta=\sigma_k^{-1}(i)-i$. If $\eta\le 0$ then the required
statement is vacuously true, and if $\eta>0$ then this implies
that swaps at positions $i,\ldots, i+\eta-1$ have appeared in this
order (with possibly other swaps in between). Thus $\eta \le
{\mathcal R_k}$, as required.
\end{proof}

\begin{proof}[Proof of \thmref{octagon}]
We use  Lemma \ref{det first row} and Corollary~\ref{C:first_row}.
With the notation there, the Edelman-Greene bijection
(\thmref{EG_bij}) shows that if $\omega=EG(T)$ then $\mathcal
R_k(\omega)=R_k(T)$.  Therefore Corollary~\ref{C:first_row} gives
for any $\eps>0$,
\[
\pusn^n \left(\forall k :\mathcal R_k < d_k+\eps n\right) \xlim 1,
\]
and by Lemma \ref{det first row} we deduce
\[
\pusn^n\left(\forall j,k : \sigma_k^{-1}(j)-j < d_k+\eps n\right) \xlim 1.
\]
Since $\sigma_k$ is a permutation this is equivalent to
\[
\pusn^n\left(\forall i,k : i-\sigma_k(i) < d_k+\eps n\right) \xlim
1.
\]
The symmetries (\ref{sym1}) and (\ref{sym2}) now imply the other
three required bounds.
\end{proof}

\begin{proof}[Proof of \thmref{Holder}]
{\it Part (i)}.  Fix $\varepsilon>0$. Consider the event
\[
E=\bigg\{ \forall i,\, \forall\, 0\le j<k\le
N\,:\,|\sigma_j^{-1}(i)- \sigma_k^{-1}(i)| \le
n\sqrt\frac{2}{N}(k-j)^{1/2} + \frac{\varepsilon n}{2}\bigg\}
\]
We will show that
\begin{equation}\label{eq:hold2}
\pusn^n E \xlim 1.
\end{equation}
This will be
enough, since the effect of linear interpolation is negligible
(that is proving \eqref{eq:hold2} for all $\varepsilon>0$ implies
the required statement for all $\varepsilon>0$). Denote $M =
\lfloor \frac{\varepsilon^2 N}{128}\rfloor$ and $K=\lfloor
\frac{128}{\varepsilon^2}\rfloor$. For each integer $0\le v \le
K$, denote the event
$$ E_v = \bigg\{
\forall i,k\,:\,|\sigma_k^{-1}(i)-\sigma_{v M}^{-1}(i)| \le
\frac{\varepsilon}{8}n + n\left(\frac{2|k-vM|}{N}\right)^{1/2}
\bigg\}. $$ By \thmref{octagon} together with stationarity,
\thmref{stat-semi}(i) we get $ \pusn (E_v) \to 1 $. Since the number
of these events is fixed, we deduce
\[
\pusn \bigg[ \bigcap_{0\le v\le K} E_v \bigg] \xrightarrow[n\to\infty]{} 1.
\]

We claim that if $\omega\in \bigcap_{0\le v\le K} E_v$ then
$\omega\in E$.  For each $0\le j<k\le  N$ consider two cases.
First, it is possible that there is some $0\le v \le K$ such that
$vM \le j<k<(v+1)M$.  In this case, $\omega\in E_v$ implies that
for all $i$
\begin{equation}\label{eq:stilltrue}
|\sigma_j^{-1}(i)-\sigma_{vM}^{-1}(i)|\le \frac{\varepsilon}{8}n+
n\left(\frac{2|j-vM|}{N}\right)^{1/2} \le \frac{\varepsilon}{8}n+
\frac{\varepsilon}{8}n = \frac{\varepsilon}{4}n,
\end{equation}
where we used the fact that $j-vM < \frac{\varepsilon^2 N}{128}$.
Similarly,
\[
|\sigma_k^{-1}(i)-\sigma_{vM}^{-1}(i)|\le \frac{\varepsilon}{4}n,
\]
and therefore
\[
|\sigma_j^{-1}(i)-\sigma_k^{-1}(i)| \le \frac{\varepsilon}{4}n+
\frac{\varepsilon}{4}n = \frac{\varepsilon}{2}n.
\]

The second possibility is that for some $0< v\le K$ we have that
$(v-1)M < j < vM \le k$. In that case, \eqref{eq:stilltrue} is
still true, and furthermore since $\omega\in E_v$ and $k-vM<k-j$
we get
\[
|\sigma_k^{-1}(i)-\sigma_{vM}^{-1}(i)|\le \frac{\varepsilon}{8}n+
n\left(\frac{2|k-vM|}{N}\right)^{1/2}\le
\frac{\varepsilon}{8}n+n\sqrt\frac{2}{N}(k-j)^{1/2}.
\]
Combining this with \eqref{eq:stilltrue} gives
\[
|\sigma_j^{-1}(i)-\sigma_k^{-1}(i)|\le \frac{\varepsilon}{2}n+
n\sqrt\frac{2}{N}(k-j)^{1/2},
\]
as claimed.
\smallskip

{\it Part (ii)}. For some fixed sequence
$\varepsilon_n\to 0$, consider the set $A_n$ of continuous functions
$T:[0,1]\to[-1,1]$ satisfying
\[
\forall t,s\in[0,1]\;:\; |T(t)-T(s)| \le \sqrt{8}|t-s|+\eps_n.
\]
By part (i) we can choose $\eps_n\to 0$ so that $\P(T_{i(n)}\in A_n) \to
1$.

Let $w(T,h)=\sup\{|T(t)-T(s)|:|t-s|\le h\}$. It follows from
$\P(T_{i(n)}\in A_n )\to 1$ that
\[
\lim_{h\to 0} \limsup_{n\to \infty} \E(w(T_{i(n)},h)\wedge 1) = 0.
\]
By \cite[Theorem 16.5]{kallenberg} we have tightness of the random sequence
$T_{i(n)}$ under this condition (note that the target space $[-1,1]$ is
compact). This establishes the existence of subsequential limits.

Now if we have a weakly convergent subsequence $T_{i(n(j))}$, then
it must have the same limit as the conditioned random variables
$\widetilde{T}^j \stackrel{d}{:=}(T_{i(n(j))} \mid T_{i(n(j))} \in
A_{n(j)})$. We may realize the sequence $\{\widetilde{T}^{j}\}$ on
the same probability space so that $\widetilde{T}^j\to T$ a.s.\
\cite[Theorem 4.30]{kallenberg}. We conclude by observing that any
limit of deterministic paths $T^j\in A_{n(j)}$ is
H\"older$(\sqrt{8},\frac{1}{2})$.
\end{proof}

\section{Great circles}
\label{S:disc}

In this section we prove \thmref{great-circle}. The idea is as follows. If
a sorting network lies close to a great circle then its trajectories are
close to sine curves up to some time change. Equivalently, it is close to a
stretchable network obtained by rotating a set of points as in the remark
in the introduction. This set of points must have roughly uniform
one-dimensional projections in all directions, so its empirical measure
must be close to $\arch_{1/2}$. Finally, since the inversion number of the
resulting configurations is close to linear in the angle of rotation, the
time change mentioned above must be linear.

Here are the details.  Denote the centre of $\S_n$ by $\bc =
\left(\frac{n+1}{2},\dots,\frac{n+1}{2}\right)$ and the radius by
$R = \sqrt{\frac{n^3-n}{12}}$. Given the circle $c_n$ we may
choose a pair of orthogonal vectors ${\bf u},{\bf v}$ of length
$R$ so that the circle has the representation $c_n =
\{c_n(\theta)\}_{\theta\in \R}$ where
\[
  c_n(\theta) = {\bf c} + {\bf u}\cos \theta + {\bf v}\sin \theta.
\]

For $k\in\{0,\dots,N\}$, define a sequence $\theta_k$ (up to
addition of multiples of $2\pi$) by
\[
\theta_k = \arg \min_\theta \| \sigma_k^{-1} - c_n(\theta) \|_\infty.
\]
Thus $c_n(\theta_k)$ is the point of $c_n$ closest in $L^\infty$ to
$\sigma_k^{-1}$. W.log.\ we may choose ${\bf u}$ so that $\theta_0=0$ (this
leaves us two possibilities for $\bv$). For other $k$, the angle $\theta_k$
is uniquely determined inductively by requiring $|\theta_{k+1}-\theta_k|
< \pi$. By symmetry, $\theta_N = (2k+1)\pi$ for some integer $k$ (and we
will see that in fact $k=0$).

Fix some $\eps>0$. The condition on $\omega_n$ implies that for
$n$ large enough (depending on $\eps$),
\begin{equation} \label{eq:sigma_approx}
  \left\| \sigma_k^{-1} - c_n(\theta_k) \right\|_\infty \le \eps n
\quad\text{for all $k$.}
\end{equation}
Since
$\|\sigma_k^{-1}-\sigma_{k+1}^{-1}\|_\infty = 1$, this implies
that
\[
  \left\| c_n(\theta_{k+1}) - c_n(\theta_k) \right\|_\infty \le 1+2\eps n.
\]
Since $R\approx n^{3/2}$, simple geometry implies that for $n$
large enough we have
\[
  |\theta_{k+1}-\theta_k|
  \le 2\arcsin\left(\frac{(1+2\eps n)\sqrt{n}}{2R}\right)
  \le 8\eps,
\]
for all $k$ (the $\sqrt n$ term comes from passing from the
$L^\infty$ norm to the $L^2$ norm). Thus $\{\theta_k\}$ does not
change too quickly. In particular, there must be some $k$ so that
either $|\theta_k-\pi/2|\le 4\eps$, or $|\theta_k+\pi/2|\le 4\eps$. We
can negate $\bv$, so w.log.\ assume the former is the case.

Considering the $i$th coordinate in \eqref{eq:sigma_approx}, one
finds that
\begin{equation} \label{eq:almost_perm}
  \Big| \sigma_k^{-1}(i) -
        \big(\tfrac{n+1}{2} + u_i \cos\theta_k + v_i \sin\theta_k \big)
  \Big| \le \eps n.
\end{equation}
We would like to show that the sorting network is approximated by
motion along the circle with constant speed, i.e.\ that $\theta_k
\approx \pi k/N$. If that were the case, part (i) of
\thmref{great-circle} would follow. As it is,
\eqref{eq:almost_perm} only implies that the paths are
approximately sine curves up to a time change. The key point here
is that the same time change applies to all particles.

Define a probability measure $\nu_n$ on $\R^2$ by
\[
  \nu_n = \frac1n \sum_{i=1}^n
  \delta\left(\tfrac2n u_i, \tfrac2n v_i \right),
\]
where $\delta(x,y)$ is the delta measure at $(x,y)$. Thus $\nu_n$ is the
empirical measure for the (rescaled) coordinates of $\bu$ and $\bv$.

\begin{lemma}
  With the above notations we have the vague convergence
  $\nu_n\Longrightarrow\arch_{1/2}$.
\end{lemma}

\begin{proof}
  We first claim that $\nu_n$ is supported inside the disc of radius 2.
  Indeed, the vector $\bc+\bu$ approximates the identity permutation, and
  so (by (\ref{eq:sigma_approx})) all entries of $\frac2n \bu$ are in
  $[-1-3\eps,1+3\eps]$. For $\bv$,
  note that there is some $k$ so that $|\theta_k-\pi/2|<4\eps$, we find that
  $\bc+\bv$ approximates $\sigma_k^{-1}$, with some additional error from
  $\bu \cos \theta_k$. Thus the coordinates of $\frac2n \bv$ are all in
  $[-1-6\eps,1+6\eps]$.

  We use the continuity theorem for the multi-dimensional characteristic
  function, (see e.g.\ \cite[Theorem~5.3]{kallenberg}). Thus it suffices to
  prove pointwise convergence of the characteristic function of $\nu_n$ to
  the characteristic function of $\arch_{1/2}$. This in turn will be
  deduced from considering the one-dimensional projections of $\nu_n$.

  More precisely, \eqref{eq:almost_perm} says that for $n$ large enough,
  for all $k$ and $i$,
  \[
  \Big| \big(\tfrac2n u_i \cos \theta_k + \tfrac2n v_i \sin \theta_k\big)
                    - \big(\tfrac2n \sigma_k^{-1}(i) -1\big) \Big| \le 3\eps.
  \]
  This states that the projection of $\nu_n$ in direction $\theta_k$ can be
  coupled to the empirical measure of a permutation scaled to $[-1,1]$ so
  that they differ by at most $2\eps$. But the scaled empirical measure of
  a permutation consists of equal point masses along an arithmetic
  progression, and does not depend on the permutation. Let
  $F(x):=0\vee\tfrac{x+1}{2}\wedge 1$ be
  the distribution function of uniform measure on $[-1,1]$. If
  $P_\theta(a,b)=a\cos \theta + b\sin\theta$ denotes projection on a line
  in direction $\theta$, then we deduce that for $n$ large enough, for any
  $x\in\R$,
  \begin{equation} \label{eq:Pthetak}
    \Big| (P_{\theta_k} \nu_n)(-\infty,x] - F(x) \Big|
    \le 3\eps + \frac1n \le 4\eps.
  \end{equation}

  For an arbitrary angle $\theta$, there is necessarily some $k$ such that
  either $|\theta_k-\theta| < 4\eps$, or the same holds for $\theta+\pi$.
  Fix $x\in[-1,1]$. Note that for two angles $\phi,\psi$, and any $z\in
  \R^2$ we have $|P_{\phi}z-P_{\psi}z|\le |z| |\phi-\psi|$. Since $\nu_n$
  is supported inside the disc of radius 2, $P_\theta \nu_n$ is close to
  $P_{\theta_k} \nu_n$, and so
  \begin{equation} \label{eq:Ptheta}
  \Big| (P_\theta \nu_n)(-\infty,x] - (P_{\theta_k} \nu_n)(-\infty,x] \Big|
  \le 8\eps.
  \end{equation}
  Combining \eqref{eq:Pthetak} and \eqref{eq:Ptheta} gives that for $n$
  large enough (depending on $\eps$), for all $\theta$ and $x\in[-1,1]$,
  \[
  \Big| (P_\theta \nu_n)(-\infty,x] - F(x) \Big| \le 12\eps.
  \]
  By monotonicity of cumulative distribution funcitons, the same bound
  holds for $1\le |x| \le 2$. Since the support of $\nu_n$ is bounded in
  the disc of radius 2, for $x$ outside $[-2,2]$ we have the stronger
  identity $(P_\theta \nu_n)(-\infty,x] = F(x)$.

  We wish to compare the characteristic functions $\phi$ and $\phi_n$ of
  $\arch_{1/2}$ and $\nu_n$ respectively. Note that for any
  $\theta$, the measure $P_\theta\arch_{1/2}$ is the uniform
  measure on
  $[-1,1]$.  We have that
  \[
  (\phi-\phi_n)(r\cos\theta,r\sin\theta)
  = \int_{\R} e^{irx} P_\theta (\nu_n-\arch_{1/2})(dx).
  \]
  Integrating by parts gives
  \begin{align*}
  \Big|(\phi-\phi_n)(r\cos\theta,t\sin\theta) \Big|
  &\le \int_{\R} \Big| ir e^{irx} \big[P_\theta(\nu_n-\arch_{1/2})\big](-\infty,x]
                       \Big| dx \\
  &\le \int_{-2}^2 2 \Big| \big[P_\theta(\nu_n-\arch_{1/2})\big](-\infty,x]
                       \Big| dx \\
  &\le 96\eps.
  \end{align*}
  Since $\eps$ can be arbitrarily small, this proves pointwise convergence
  of the characteristic functions, and therefore convergence of $\nu_n$ to
  $\arch_{1/2}$.
\end{proof}

\begin{lemma} \label{L:theta_linear}
  With the above notation,
  \[
  \max_k \left| \theta_k-\frac{k\pi}{N} \right| \xrightarrow[n\to\infty]{} 0
  \]
\end{lemma}

\begin{proof}
  Fix some $\theta$, and let $P_\theta$ denote the projection on a line in
  direction $\theta$, so that
  \[
  P_\theta(x,y)=x\cos \theta+y\sin \theta,
  \]
  and consider the permutation $\rho_n(\theta)$ derived from $\bu,\bv$
  by arranging $i\in[1,n]$ in increasing order of $P_\theta(u_i,v_i)$. We
  first estimate the inversion number $\inv(\rho_n(\theta))$. Define
  \[
   A(\theta) = \left\{((x,y),(x',y')) \in(\R^2)^2 \ : \ \begin{aligned}
       x&<x', \\ P_\theta(x,y) &>P_\theta(x',y')
     \end{aligned} \right\}.
  \]
  We have
  \[
  \frac1N \inv(\rho_n(\theta)) = \iint \ind_{A(\theta)}\; d\nu_n\; d\nu_n,
  \]
  which is a continuous functional of $\nu_n$. Since $\nu_n\Longrightarrow
  \arch_{1/2}$, this implies
  \begin{equation}\label{This}
    \frac1N \inv(\rho_n(\theta)) \xrightarrow[n\to\infty]{}
    \iint \ind_{A(\theta)} \;d\arch_{1/2}\; d\arch_{1/2} =
    \frac{\theta}{\pi}
  \end{equation}
  To check the last equality, note that the integral is the probability
  that the $x$-projections of two points chosen independently from
  $\arch_{1/2}$ change order after rotation by at most $\theta$. By
  rotational invariance the angle of the line between two such points is
  uniform on $[0,\pi]$.

  Equation (\ref{This}) holds for any fixed $\theta$. However,
  $\inv(\rho_n(\theta))$ is increasing in $\theta\in[0,\pi]$, and
  consequently,
  \begin{equation}\label{eq:inv_rho}
    \frac1N \inv(\rho_n(\theta)) \xrightarrow[n\to\infty]{}
    \frac\theta \pi \quad\text{uniformly in $\theta\in[0,\pi]$}.
  \end{equation}

  Comparing \eqref{eq:almost_perm} for $i,j$ we find that for any $\eps>0$,
  for $n$ large enough we have
  \begin{equation}
    \label{eq:ij_order}
  \sigma_k^{-1}(i)-\sigma_k^{-1}(j)
  = (u_i-u_j)\cos\theta_k + (v_i-v_j)\sin\theta_k + \delta
  \end{equation}
  with $|\delta|<2\eps n$ for all $i,j$.
  Say a pair $(i,j)$ is $\theta${\bf-uncertain} if
  \[
    |(u_i-u_j)\cos\theta_k + (v_i-v_j)\sin\theta_k| \le 2\eps.
  \]

  So, for any $i<j$ we have that $\sigma_k^{-1}(i)>\sigma_k^{-1}(j)$ if and
  only if $(\rho_n(\theta_k))(i)>(\rho_n(\theta_k))(j)$, unless $(i,j)$ is
  $\theta_k$-uncertain. Recall that any such pair $i<j$ with
  $\sigma_k^{-1}(i)>\sigma_k^{-1}(j)$ contributes 1 to the number of
  inversions of $\sigma_k^{-1}$ (hence, if $(i,j)$ is not
  $\theta_k$-uncertain, also to $\text{inv}(\rho_n(\theta_k))$).
  Consequently, $\inv(\rho_n(\theta_k))$ differs from
  $\inv(\sigma_k^{-1})=k$ by at most the number of $\theta_k$-uncertain
  pairs.

  It remains to bound the number of $\theta$-uncertain pairs. Fix $\eps>0$,
  and consider the set $S_{4\eps}$ of all strips of width $4\eps$ in
  $\R^2$. Since $\nu_n\Longrightarrow \arch_{1/2}$, we have
  \[
  \limsup_{n\to\infty} \sup_{A\in S_{4\eps}} \nu_n(A) \le 2\eps
  \]
  because $\arch_{1/2}(A)\le 2\eps$ for any such strip. This implies
  that for large $n$ for any $i$ there are at most $2\eps n$ values of
  $j$ such that $(i,j)$ is $\theta$-uncertain for some $\theta$. In
  summary, for $n$ large enough, depending only on $\eps$, and any
  $\theta$, the total number of $\theta$-uncertain points is at most $2\eps
  n^2$.

  Combining \eqref{eq:ij_order} and the above discussion we find that for
  large $n$
  \[
  |\inv(\sigma_k)-\inv(\rho_n(\theta_k))| \le 2\eps n^2,
  \]
  uniformly in $\theta$. Since $\inv(\sigma_k)=k$, combining with
  \eqref{eq:inv_rho} yields the result.
\end{proof}

\begin{proof}[Proof of \thmref{great-circle}]
  Combining \lemref{theta_linear} and \eqref{eq:almost_perm} gives part
  (i):
  \begin{equation} \label{eq:uv_sin}
  \max_{i,t} \left| (\tfrac 2n \sigma_{[tN]}(i) -1)
                    - \tfrac2n(u_i\cos(\pi t)+v_i\sin (\pi t)) \right|
                  \xrightarrow[n\to\infty]{} 0.
  \end{equation}
  In particular, we have $2u_i/n=2i/n-1+o(1)$. By inserting
  \eqref{eq:uv_sin} into the definition of $\mu_t$ in \eqref{eq:mutdef}, we
  find that $\mu_t(\omega_n)$ is close to $R_t \nu_n$, where $R_t$ is the
  linear map $R_t(x,y) = (x,x\cos(\pi t) + y\sin(\pi t))$, in the sense the
  two measures can be coupled with maximal distance tending to 0. Since
  $\nu_n \Longrightarrow \arch_{1/2}$, this implies $\mu_t(\omega_n)
  \Longrightarrow \arch_t$, which is (ii).

  Next, we prove (iii). To sample from the scaled swap process
  $\eta(\omega_n)$ one may choose uniformly a pair of particles $i,j$ and
  consider the time and location of their swap. Consider the pair of points
  $z_i=\frac2n(u_i,v_i)$ and $z_j=\frac2n(u_j,v_j)$. If $i,j$ are swapped
  at step $k$ of the network, then $\sigma_k^{-1}(i)-\sigma_k^{-1}(j)=1$,
  so by (\ref{eq:ij_order}) we have for $n$ large enough that
  $|P_{\theta_k}(z_i)-P_{\theta_k}(z_j)|\le 2\varepsilon$, and by
  (\ref{eq:almost_perm}) the scaled location of the swap is given to within
  $2\varepsilon$ by $P_{\theta_k} z_i$. Thus for any $i,j$, the time of the
  $(i,j)$ swap is given by the angle of a certain line, and the location of
  the swap by the distance of the line from the origin, where this line
  passes within distance $\eps$ of both $z_i$ and $z_j$. Thus, unless
  $z_i,z_j$ are sufficiently close, the location of the two points
  approximately determines the time and place of the swap.

  Specifically, for any pair $z,z'$, as $\eps\to 0$ the set of possible
  times converges to a single time, and the set of possible locations
  converges to a single location. Since $\nu_n\Longrightarrow\arch_{1/2}$,
   it follows that
  $\eta(\omega_n)$ converges to the measure resulting from applying the
  same operation to $\arch_{1/2}$.

  Let $z,z'\in\R^2$ be chosen independently with law $\arch_{1/2}$. Let
  $\theta\in[0,\pi]$ be the angle that the line through them makes with the
  positive $y$-axis, and let $r:=z_1\cos\theta+z_2\sin\theta$ be its signed
  distance from the origin. It remains to prove that $\theta$ and $r$ are
  independent, $\theta$ is uniform in $[0,\pi]$ and $r$ has law $\semi$.
  Independence and uniformity of $\theta$ are clear by rotational symmetry
  of $\arch_{1/2}$. Finally, to calculate the distribution of $r$ we
  introduce some further variables. Let $\widehat z,\widehat z'$ be $z,z'$
  rotated by $-\theta$ and let $\widehat z=(r,y)$ and $\widehat z'=(r,y')$
  be their coordinates. Let $w=y/\sqrt{1-r^2}$ and $w'=y'/\sqrt{1-r^2}$.
  Thus we have
  \begin{align*}
    z_1&=r\cos\theta-w\sqrt{1-r^2}\sin\theta \\
    z_2&=r\sin\theta+w\sqrt{1-r^2}\cos\theta \\
    z_1'&=r\cos\theta-w'\sqrt{1-r^2}\sin\theta \\
    z_2'&=r\sin\theta+w'\sqrt{1-r^2}\cos\theta.
  \end{align*}
  We can compute the probability density function of $r$ using the Jacobian
  of the transformation $(z_1,z_2,z_1',z_2')\mapsto(r,\theta,w,w')$; after
  some straightforward manipulation we obtain
  \begin{gather*}
    \int_{0}^\pi \int_{-1}^{1} \int_{-1}^{1} \left|
    \frac{\partial(z_1,z_2,z_1',z_2')}{\partial(r,\theta,w,w')}
    \right|
    \frac{1}{2\pi\sqrt{1-\|z\|_2^2}}\;\frac{1}{2\pi\sqrt{1-\|z'\|_2^2}}
    \;dw\,dw'\,d\theta \\
    =\tfrac{2}{\pi}\sqrt{1-r^2}
  \end{gather*}
  as required.
\end{proof}

\subsection*{Acknowledgements}

We thank Nathana\"el Berestycki, Alex Gamburd, Alan Hammond, Pawel
Hitczenko, Martin Kassabov, Rick Kenyon, Scott Sheffield, David Wilson
and Doron Zeilberger for many valuable conversations.  This work has
benefitted greatly from the superb resources provided by a meeting at
BIRS (Banff, Canada) and by the 2005 programme in Probability,
Algorithms and Statistical Physics at MSRI (Berkeley, USA).

\bibliographystyle{abbrv}
\bibliography{sort}

\def\cprime{$'$}
\begin{thebibliography}{10}

\bibitem{sort-pen}
O.~Angel, A.~E. Holroyd, and D.~Romik.
\newblock Directed random walk on the permutahedron.
\newblock In preparation.

\bibitem{sort-local}
O.~Angel, A.~E. Holroyd, and B.~Virag.
\newblock The local limit of the uniform sorting network.
\newblock In preparation.

\bibitem{sort-weak}
O.~Angel, M.~Kassabov, A.~E. Holroyd, D.~Romik, and B.~Virag.
\newblock Bounds on the uniform sorting network.
\newblock In preparation.

\bibitem{Biane}
P.~Biane.
\newblock Representations of symmetric groups and free probability.
\newblock {\em Adv. Math.}, 138(1):126--181, 1998.

\bibitem{edelmangreene}
P.~Edelman and C.~Greene.
\newblock Balanced tableaux.
\newblock {\em Adv. in Math.}, 63(1):42--99, 1987.

\bibitem{feller}
W.~Feller.
\newblock {\em An introduction to probability theory and its applications.
  {V}ol. {I}}.
\newblock Third edition. John Wiley \& Sons Inc., New York, 1968.

\bibitem{felsner}
S.~Felsner.
\newblock The skeleton of a reduced word and a correspondence of {E}delman and
  {G}reene.
\newblock {\em Electron. J. Combin.}, 8(1):Research Paper 10, 21 pp.
  (electronic), 2001.

\bibitem{finch}
S.~R. Finch.
\newblock {\em Mathematical constants}, volume~94 of {\em Encyclopedia of
  Mathematics and its Applications}.
\newblock Cambridge University Press, Cambridge, 2003.

\bibitem{framerobinsonthrall}
J.~S. Frame, G.~B. Robinson, and R.~M. Thrall.
\newblock The hook graphs of the symmetric groups.
\newblock {\em Canadian J. Math.}, 6:316--324, 1954.

\bibitem{garsia}
A.~Garsia.
\newblock The saga of reduced factorizations of elements of the symmetric
  group.
\newblock Preprint, \\
  \verb+http://www.math.ucsd.edu/~garsia/recentpapers/saga.pdf+.

\bibitem{goodman}
J.~E. Goodman and J.~O'Rourke, editors.
\newblock {\em Handbook of discrete and computational geometry}.
\newblock Discrete Mathematics and its Applications (Boca Raton). Chapman \&
  Hall/CRC, Boca Raton, FL, second edition, 2004.

\bibitem{greeneetal}
C.~Greene, A.~Nijenhuis, and H.~S. Wilf.
\newblock A probabilistic proof of a formula for the number of {Y}oung tableaux
  of a given shape.
\newblock {\em Adv. in Math.}, 31(1):104--109, 1979.

\bibitem{kallenberg}
O.~Kallenberg.
\newblock {\em Foundations of modern probability}.
\newblock Probability and its Applications (New York). Springer-Verlag, New
  York, second edition, 2002.

\bibitem{Ke1}
S.~Kerov.
\newblock A differential model for the growth of {Y}oung diagrams.
\newblock In {\em Proceedings of the St. Petersburg Mathematical Society, Vol.
  IV}, volume 188 of {\em Amer. Math. Soc. Transl. Ser. 2}, pages 111--130,
  Providence, RI, 1999. Amer. Math. Soc.

\bibitem{Ke2}
S.~V. Kerov.
\newblock Transition probabilities of continual {Y}oung diagrams and the
  {M}arkov moment problem.
\newblock {\em Funktsional. Anal. i Prilozhen.}, 27(2):32--49, 96, 1993.

\bibitem{knuth}
D.~E. Knuth.
\newblock {\em The Art of Computer Programming, Vol. 3: Sorting and Searching}.
\newblock Addison-Wesley Publishing Co., Reading, Mass.-London-Don Mills, Ont.,
  1973.
\newblock Addison-Wesley Series in Computer Science and Information Processing.

\bibitem{lascouxschutz}
A.~Lascoux and M.-P. Sch{\"u}tzenberger.
\newblock Structure de {H}opf de l'anneau de cohomologie et de l'anneau de
  {G}rothendieck d'une vari\'et\'e de drapeaux.
\newblock {\em C. R. Acad. Sci. Paris S\'er. I Math.}, 295(11):629--633, 1982.

\bibitem{little}
D.~P. Little.
\newblock Combinatorial aspects of the {L}ascoux-{S}ch\"utzenberger tree.
\newblock {\em Adv. Math.}, 174(2):236--253, 2003.

\bibitem{pittelromik}
B.~Pittel and D.~Romik.
\newblock Limit shapes for random square young tableaux.
\newblock {\em Adv. Appl. Math.}, to appear.

\bibitem{reiner}
V.~Reiner.
\newblock Note on the expected number of {Y}ang-{B}axter moves applicable to
  reduced decompositions.
\newblock {\em European J. Combin.}, 26(6):1019--1021, 2005.

\bibitem{stanleypaper}
R.~P. Stanley.
\newblock On the number of reduced decompositions of elements of {C}oxeter
  groups.
\newblock {\em European J. Combin.}, 5(4):359--372, 1984.

\bibitem{stanley}
R.~P. Stanley.
\newblock {\em Enumerative combinatorics. {V}ol. 2}, volume~62 of {\em
  Cambridge Studies in Advanced Mathematics}.
\newblock Cambridge University Press, Cambridge, 1999.
\newblock With a foreword by Gian-Carlo Rota and appendix 1 by Sergey Fomin.

\bibitem{white}
A.~T. White.
\newblock Fabian {S}tedman: the first group theorist?
\newblock {\em Amer. Math. Monthly}, 103(9):771--778, 1996.

\end{thebibliography}

\vspace{6mm}
\begin{minipage}{12cm}

Omer Angel: {\tt angel@utstat.toronto.edu} \\
Department of Statistics, University of Toronto,\\
100 St George St., Toronto ON M5S 3G3, Canada.

\vspace{3mm} \noindent Alexander E. Holroyd: {\tt holroyd@math.ubc.ca} \\
Department of Mathematics, University of British Columbia, \\
121-1984 Mathematics Rd., Vancouver BC V6T 1Z2, Canada.

\vspace{3mm} \noindent Dan Romik: {\tt romik@stat.berkeley.edu}\\
Department of Statistics, University of California, \\
367 Evans Hall, Berkeley CA 94720-3860, USA.

\vspace{3mm} \noindent B\'alint Vir\'ag: {\tt balint@math.toronto.edu}\\
Department of Mathematics, University of Toronto,\\
40 St George St., Toronto ON M5S 2E4, Canada.
\end{minipage}

\end{document}